\documentclass{article}

\usepackage{float}
\usepackage{cite}
\usepackage[utf8x]{inputenc}
\usepackage{amsmath}
\usepackage{amssymb}
\usepackage{amsthm}
\usepackage{graphicx}
\usepackage[colorinlistoftodos]{todonotes}
\usepackage[colorlinks=true, allcolors=blue]{hyperref}
\usepackage{caption}
\usepackage{subcaption}
\usepackage{authblk}
\usepackage{bbm}
\usepackage{booktabs}
\usepackage{multirow}
\usepackage{array}
\usepackage[left=4cm,right=4cm,top=2.5cm,bottom=2cm,includeheadfoot]{geometry}
\usepackage{mathtools}
\usepackage{siunitx} 
\DeclareSIUnit\Molar{\text{M}}
\DeclarePairedDelimiter\abs{\lvert}{\rvert}

\newcommand{\cs}[1]{\textcolor{red}{cs: #1}}

\theoremstyle{definition}

\title{Rate-limiting recovery processes in neurotransmission under sustained stimulation}

\author[a]{Ariane Ernst}
\author[a]{Nathalie Unger}
\author[a,b]{Christof Sch\" utte}
\author[c,d,+]{Alexander Walter\footnote{awalter@sund.ku.dk}}
\author[a,+]{Stefanie Winkelmann\footnote{winkelmann@zib.de}}

\affil[a]{Zuse Institute Berlin, Berlin, Germany}
\affil[b]{Freie Universit\" at Berlin, Faculty of Mathematics and Computer Science, Berlin, Germany}
\affil[c]{University of Copenhagen, Faculty of Health and Medical Sciences, Department of Neuroscience, Copenhagen, Denmark}
\affil[d]{Leibniz-Forschungsinstitut f\"ur Molekulare Pharmakologie, FMP im Charit\'eCrossOver, Berlin, Germany}
\affil[+]{these authors contributed equally to this work}

\begin{document}

\maketitle

\begin{abstract} 
At chemical synapses, an arriving electric signal induces the fusion of vesicles with the presynaptic membrane, thereby releasing neurotransmitters into the synaptic cleft. 
After a fusion event, both the release site and the vesicle undergo a recovery process before becoming available for reuse again. Of central interest is the question which of the two restoration steps acts as the limiting factor during neurotransmission under high-frequency sustained stimulation. 
In order to investigate this question, we introduce a novel non-linear reaction network which 
involves explicit recovery steps for both the vesicles and the release sites, and includes the induced time-dependent output current.
The associated reaction dynamics  
are formulated by means of ordinary differential equations (ODEs), as well as via the associated stochastic jump process.
While the stochastic jump model describes a single release site, the average over many release sites is close to the ODE solution and shares its periodic structure. The reason for this can be traced back to the insight that recovery dynamics of vesicles and release sites are statistically almost independent. 
A sensitivity analysis on the recovery rates based on the ODE formulation reveals that neither the vesicle nor the release site recovery step can be identified as the essential rate-limiting step but that the rate-limiting feature changes
over the course of stimulation. 
Under sustained stimulation the dynamics given by the ODEs exhibits transient changes leading from an initial depression of the postsynaptic response to an asymptotic periodic orbit, while the individual trajectories of the stochastic jump model lack the oscillatory behavior and asymptotic periodicity of the ODE-solution. 
 

\quad 

\textbf{Key words:} nonlinear reaction networks, neurotransmission models, vesicle fusion dynamics, sustained stimulation
\end{abstract}






\section{Introduction}\label{sec:mot}
Communication in the nervous system relies on chemical transmission across synapses. For this, neurotransmitters are released from presynaptic neurons by the fusion of transmitter-containing synaptic vesicles with the plasma membrane. The liberated transmitter is detected by postsynaptic receptors which induces a response. At the presynapse, evoked transmitter release is limited to so-called release sites at which synaptic vesicles attach to the plasma membrane (a process referred to as \textit{vesicle docking}) and functionally mature to become responsive to presynaptic stimulation (a process referred to as \textit{vesicle priming}) \cite{sudhof2004synaptic,kaeser2017readily,verhage2008,walter_vesicle_2018}. Transmitter release is typically induced by presynaptic action potentials, i.e., brief de- and re-polarisations of the cellular membrane potential that lead to the opening of voltage gated $\text{Ca}^{2+}$ ion channels \cite{sudhof2013neurotransmitter,catterall2011voltage,stanley2016nanophysiology}. 
The resulting elevation of the presynaptic $\text{Ca}^{2+}$ concentration following $\text{Ca}^{2+}$ influx through these channels triggers synaptic vesicle fusion by activating the vesicular $\text{Ca}^{2+}$ sensing protein Synaptotagmin \cite{sudhof2013neurotransmitter,koh2003synaptotagmin,kobbersmed2022allosteric}. 

During high-frequency sustained stimulation, most synapses exhibit a depression of the initial postsynaptic response to a plateau \cite{zucker2002short,kavalali2006synaptic,gardner1980rate,galarreta1998frequency}. 
Continued activity puts a great demand on the cell to replenish the active zone in time with 
synaptic vesicles as well as release site proteins.  Both of these are finite resources that may be expended  during continued exocytosis and therefore need to be replenished for sustained activity. The observed depression of the postsynaptic current during prolonged stimulation may thus be explained by the refractory recycling of the release sites and/or the depletion of available, fusion-competent synaptic vesicles and by the time it takes to replenish those. 

The membrane of synaptic vesicles that underwent fusion is taken up by endocytosis, and further processing and neurotransmitter re-uptake is required to regenerate a new synaptic vesicle \cite{kononenko2015molecular}. 
 Initially, it was thought that endocytosis itself was slow (on the timescale of  tens of seconds  \cite{granseth2006clathrin}), but more recently it became clear that endocytosis can occur much faster ("ultrafast" millisecond timescale) at presynaptic membranes \cite{watanabe2013ultrafast,delvendahl2016fast}. However, the full process of synaptic vesicle reformation is  thought to take longer  \cite{lenzi2002depolarization,yamashita2018vesicular,nakakubo2020vesicular,saheki2012synaptic}, which is why vesicle replenishment is often assumed to be the limiting step during sustained stimulation \cite{kavalali2006synaptic,sudhof2004synaptic,katz1993neurotransmitter,betz1970depression,betz1993optical, wu1999reduced, fernandez2004kinetics, rizzoli2005synaptic}.
On the other hand, there is usually a large supply of reserve vesicles in a synaptic bouton and rapid replenishment from a large reserve pool could counteract synaptic depression (or even cause facilitation), and recent experimental data indicate that vesicle replenishment may commence much faster than previously thought, within milliseconds \cite{miki2016actin,sudhof2004synaptic,neher2010rate,neher2017some,kaeser2017readily,magleby1987short}. 

Apart from vesicle resupply, the release sites themselves may set constraints on further presynaptic activity, for instance, if they need to undergo some form of clearance and/or recycling before taking up another vesicle. Experiments in Drosophila, where the endocytosis machinery was acutely blocked, demonstrated that impairing endocytosis affected repetitive synaptic activity on a timescale of milliseconds, indicating that fast site clearance by endocytosis could be a major factor to maintain synaptic activity \cite{kawasaki_fast_2000}. Accordingly, release site recycling was  estimated to happen on very short timescales, within tenths of a second \cite{neher2010rate}. 
It is currently not known which of the two reactions -- vesicle replenishment or release site resupply -- is limiting sustained synaptic activity (as pointed out earlier\cite{neher2010rate}). Experimentally this is difficult to distinguish as most read-outs with sufficient temporal resolution (e.g. electrophysiology, live imaging) quantify the downstream neurotransmitter release  and cannot directly report on upstream processes. More recently, rapid, high pressure freeze fixation shortly after synaptic stimulation provided insight into the morphological changes and provided a first account of the kinetics of vesicle reformation at neurotransmitter release sites \cite{watanabe2013ultrafast}. Yet, even such approaches cannot resolve whether this reformation is limited by the vesicular association to the sites or the availability of the sites to receive a vesicle. Thus, to date it is not known to which degree either (or both) of the reactions limit continued synaptic output, or whether this can even be distinguished. An insight into this would be valuable, for example to understand which effects to expect if pathogenic or environmental factors selectively affect them.

\quad

In this paper, we set out to investigate to which extent either vesicle replenishment or release site recycling limit neural activity during sustained activation. Based on the unpriming model investigated in prior work \cite{kobbersmed2020rapid,ernst2022variance}, we introduce a novel non-linear model that includes the combined recycling dynamics of vesicles and release sites.
Given the underlying reaction network, we first describe the dynamics  by a set of ordinary differential equations (ODEs) and include the
postsynaptic output by convolving the fusion events with a characteristic postsynaptic response evoked by a single vesicle \cite{kobbersmed2020rapid}. 
The parameter values have been estimated based on literature with the aim to let the ODE solution describe the average postsynaptic response signal at the Drosophila neuromuscular junction. 
Typical solutions of the ODE model under sustained stimulation exhibit transient dynamics leading from an initial depression of the postsynaptic response to an asymptotic periodic orbit. We demonstrate that the asymptotic periodic solution oscillates around a unique steady state given by the running average of parameters under sustained stimulation.


The ODE model with these parameters is then used to investigate the impact of the two recovery rates (vesicles vs release sites) on the dynamics.
As a measure for the influence of both recovery steps we choose the sensitivity of the output current with respect to the recycling rates. 
We show that the identity of the rate-limiting process depends on the point in time during stimulation:
With the investigated parameter values, the neural output during 100 Hz stimulation is initially more sensitive to the release site replenishment. Later (once vesicles have accumulated in the recycling state) this shifts to a high sensitivity to the vesicle replenishment. We observe that this behaviour is conserved over a large range of parameters. 

Next, we extend our analysis to the stochastic jump process model given by the associated reaction network. We observe that the characteristic transient behavior and asymptotic periodicity of the ODE dynamics under sustained stimulation is not visible for an individual stochastic trajectory. 
This is no surprise since the stochastic jump process is describing a single release site and its discrete stochastic response to stimulation. 
However, the transient dynamics and asymptotic periodicity return when considering the junction current averaged over many release sites. This averaged current
converges to the first-order moment of the stochastic output current, which is demonstrated to be in very close agreement with the ODE-solution. We trace this similarity back to the statistical independence of the recovery processes
and a resulting small correlation between release site and vesicle supply. The
agreement is independent of the model parameters and therefore supports the validity
and applicability of the ODE-based sensitivity analysis.

\quad

Our paper is organized as follows. In Sec.~\ref{sec:ODE} we introduce the recovery model as a reaction network including explicit recovery steps for vesicles and release sites. Next we numerically analyze the system response to sustained stimulation, including sensitivity analysis of the two recovery processes. In Sec.~\ref{sec:stochastic} we extend our analysis to stochastic dynamics. The total junction current induced by several release sites is simulated and compared with the ODE-solution. Finally, the system's first- and second-order moments and its correlation function are investigated.

\section{Recovery dynamics of vesicles and release sites}
\label{sec:ODE}

In this section, we introduce the recovery model for the interaction dynamics of vesicles and release sites. The dynamics are formulated in terms of an ODE (more precisely, the reaction rate equation), which is solved numerically in order to study the temporal evolution of the system's response to sustained stimulation. By a sensitivity analysis, we investigate the impact of the recovery rates onto the output current. 

\subsection{ODE model for many release sites} \label{sec:simplemodel}

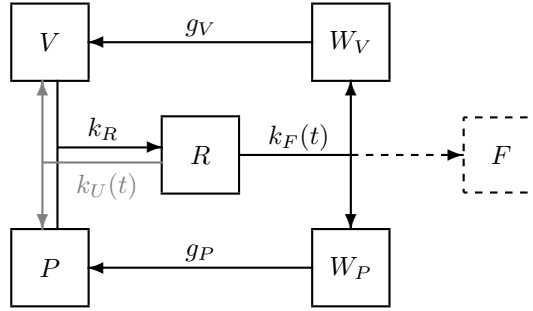
\begin{figure}
\begin{center}
\begin{picture}(6,4)
\thicklines
\put(0,0){\framebox(1,1){$P$}}
\put(4,0){\framebox(1,1){$W_P$}}
\put(0,3){\framebox(1,1){$V$}}
\put(4,3){\framebox(1,1){$W_V$}}
\put(2,1.5){\framebox(1,1){$R$}}
\put(6,1.5){\dashbox{0.1}(1,1){$F$}}
\put(4,0.5){\vector(-1,0){3}}
\put(4,3.5){\vector(-1,0){3}}
\put(0.6,2.1){\vector(1,0){1.4}}
\put(0.6,1){\line(0,1){2}}
\put(0.4,1.9){{\color{gray}\line(1,0){1.6}}}
\put(0.4,2){{\color{gray}\vector(0,1){1}}}
\put(0.4,2){{\color{gray}\vector(0,-1){1}}}
\put(4.5,2){\vector(0,1){1}}
\put(4.5,2){\vector(0,-1){1}}
\put(3,2){\line(1,0){1.5}}
\multiput(4.5,2)(0.2,0){7}{\line(1,0){0.1}} 
\put(5.9,2){\vector(1,0){0.1}}
\put(2.3,0.65){$g_P$}
\put(2.3,3.65){$g_V$}
\put(3.4,2.15){$k_F(t)$}
\put(1,2.25){$k_R$}
\put(0.85,1.5){{\color{gray}$k_U(t)$}}
\end{picture}
\end{center}
\caption{\textbf{The recovery model.} A freely available vesicle $V$ may attach to (and detach from) a freely available release site $P$.  From the resulting joint state $R$, fusion can occur, transferring both the vesicle and the release site into their recovery states $W_V, W_P$, respectively, and increasing the number $F$ of fusion events. By recovery, the vesicles and release sites turn back into the states $V$ and $P$, respectively.}
\label{fig:simplified_model}
\end{figure}

Based on the unpriming model introduced by Kobbersmed et al.~\cite{kobbersmed2020rapid}, see Sec.~\ref{parameter_estimation} for a short summary, we introduce the following recovery model for the combined recycling dynamics of a large number of vesicles and release sites, see Fig.~\ref{fig:simplified_model} for an illustration of the underlying reaction network. 
Note that experimentally measured currents are also the results of the combined activity of numerous release sites. 
Later, in Sec.~\ref{sec:stochastic} below, we will consider individual release sites and discuss the differences and similarities between the two cases.

In the model, each release site can be in three different states: It can be freely available (state $P$), or there is a vesicle attached to it (both together forming the complex $R$), or the release site is in a recovery state $W_P$. Similarly, there are three states for each vesicle: It can be freely available (state $V$), or attached to a release site (joint state $R$), or in recovery (state $W_V$). A freely available vesicle can dock with a certain rate $k_R>0$ to a freely available release site, which is expressed by the second-order reaction 
\begin{equation}\label{reaction1}
    V+P \stackrel{k_R}{\longrightarrow} R.
\end{equation}
This reaction is reversible by an umpriming reaction of the form
\begin{equation}\label{reaction2}
    R \stackrel{k_U(t)}{\longrightarrow} V+P,
\end{equation}
meaning that the vesicle detaches from the release site again. This happens at a time-dependent rate $k_U(t)\geq 0$. 
The docked vesicle may fuse with the membrane, thereby transferring both itself and the release site into the recovery state,
\begin{equation}\label{reaction3}
    R \stackrel{k_F(t)}{\longrightarrow} W_V + W_P,
\end{equation}
for a time-dependent fusion rate $k_F(t)\geq 0$.
Independently of each other, the vesicle and the release site recover according to the reactions
\begin{equation}\label{reaction4}
    W_V \stackrel{g_V}{\longrightarrow} V, \quad W_P \stackrel{g_P}{\longrightarrow} P,
\end{equation}
for time-independent rates $g_V,g_P>0$, respectively. 

The cumulative state of the system at time $t\geq 0$ is given by 
\begin{equation}\label{X}
    \boldsymbol{X}(t)= (X_i(t))_{i=1,...,5} = \left(V(t),W_V(t),W_P(t),R(t),P(t)\right)^\top
\end{equation}
where $X_i(t)$ stands for the (average) number of vesicles or release sites in the respective state. Additionally, there is the counting process $(F(t))_{t\geq 0}$ with $F(t)$ referring to the number of fusion events \eqref{reaction3} up to time $t$.


The dynamics are described by the reaction rate equation 
\begin{equation}\label{RRE}
    \Dot{\boldsymbol{X}}(t)=h(\boldsymbol{X}(t),t)
\end{equation}
with 
\begin{equation}\label{RRE:simple}
h(\boldsymbol{X}(t),t):=\left(\begin{matrix}
-k_R V(t) P(t) +g_V W_V(t) +k_U(t)R(t)\\
k_F(t) R(t) -g_V W_V(t)\\
k_F(t) R(t) - g_P W_P(t)\\
k_R V(t) P(t)  -  k_F(t) R(t) -k_U(t)R(t)\\
- k_R V(t) P(t) +g_P W_P(t) +k_U(t)R(t)
\end{matrix}\right).
\end{equation}
It is straightforward to see that both the total number of release sites and the total number of vesicles are conserved in the course of time, i.e., given that the initial states fulfill $ R(0)+P(0)+W_P(0) = n_{\text{sites}} $ and $ R(t)+V(t)+W_V(t) =  n_{\text{ves}}$ for  $n_{\text{sites}},n_{\text{ves}} \in \mathbb{N}_+$, we have two conservation laws 
\begin{equation}\label{n_constant}
    R(t)+P(t)+W_P(t) =  n_{\text{sites}}, \quad  R(t)+V(t)+W_V(t) =   n_{\text{ves}}
\end{equation}
for all $t\geq 0$, making the system effectively three-dimensional. 
The number $F(t)$ of fusion events is set to fulfill $F(0)=0$ and
\begin{equation}\label{eq:dFdt}
    \frac{d}{dt}F(t) = k_F(t)R(t).
\end{equation}
The postsynaptic \textit{response current} (output signal) induced by the dynamics of the process $(\boldsymbol{X}(t))_{t\geq 0}$ is given by the convolution of the derivative $\dot F(t) = \frac{d}{dt}F(t)$ of $F(t)$ with an impulse response function $g$  \cite{kobbersmed2020rapid,ernst2022variance}:
\begin{equation}\label{output_signal}
    C(t) := (\dot F \ast g)(t) = \int_{-\infty}^\infty \dot F(s)g(t-s) ds,
\end{equation}
where $g$ is specified by Eq.~\eqref{eq:g} in the Appendix.



\paragraph{Symmetry of the model.}
At first sight, the recovery model seems to be symmetric in the sense that the recovery dynamics of the release sites have the same structure as those of the vesicles. It is not directly evident why the roles of release sites and vesicles should not be interchangeable. However, this seemingly symmetric situation is broken by the fact that the numbers of vesicles and release sites are different and the values of the corresponding recovery rates differ, too: Per release site there are typically several vesicles, each of which needs more time for recovery than the release site itself, see the parameter estimation in Sec.~\ref{parameter_estimation}. 

\paragraph{Steady state. }
When assuming all rates to be time-independent (especially $k_U(t)=\text{const}$ and $k_F(t)=\text{const}$), the function $h$ defined in \eqref{RRE:simple} also becomes explicitly independent of time, and 
the process $\boldsymbol{X}(t)$ given by the reaction rate equation $\Dot{\boldsymbol{X}}(t)=h(\boldsymbol{X}(t))$ has a unique steady state which continuously depends on the values of the reaction rates and the numbers $n_{\text{sites}}$ and $n_{\text{ves}}$, as demonstrated in Sec.~\ref{app:steady_state}. I.e., there is exactly one state $\hat{\boldsymbol{x}} \in \mathbb{R}_+^5$ with $h(\hat{\boldsymbol{x}})=0$, and the process will approach this state asymptotically, no matter where it starts from (as long as the initial state is non-negative). 


\subsection{System response to sustained stimulation} \label{sec:sustained_stimulation}

\begin{figure}[ht!]
  \centering
   {\includegraphics[width=0.8\textwidth]{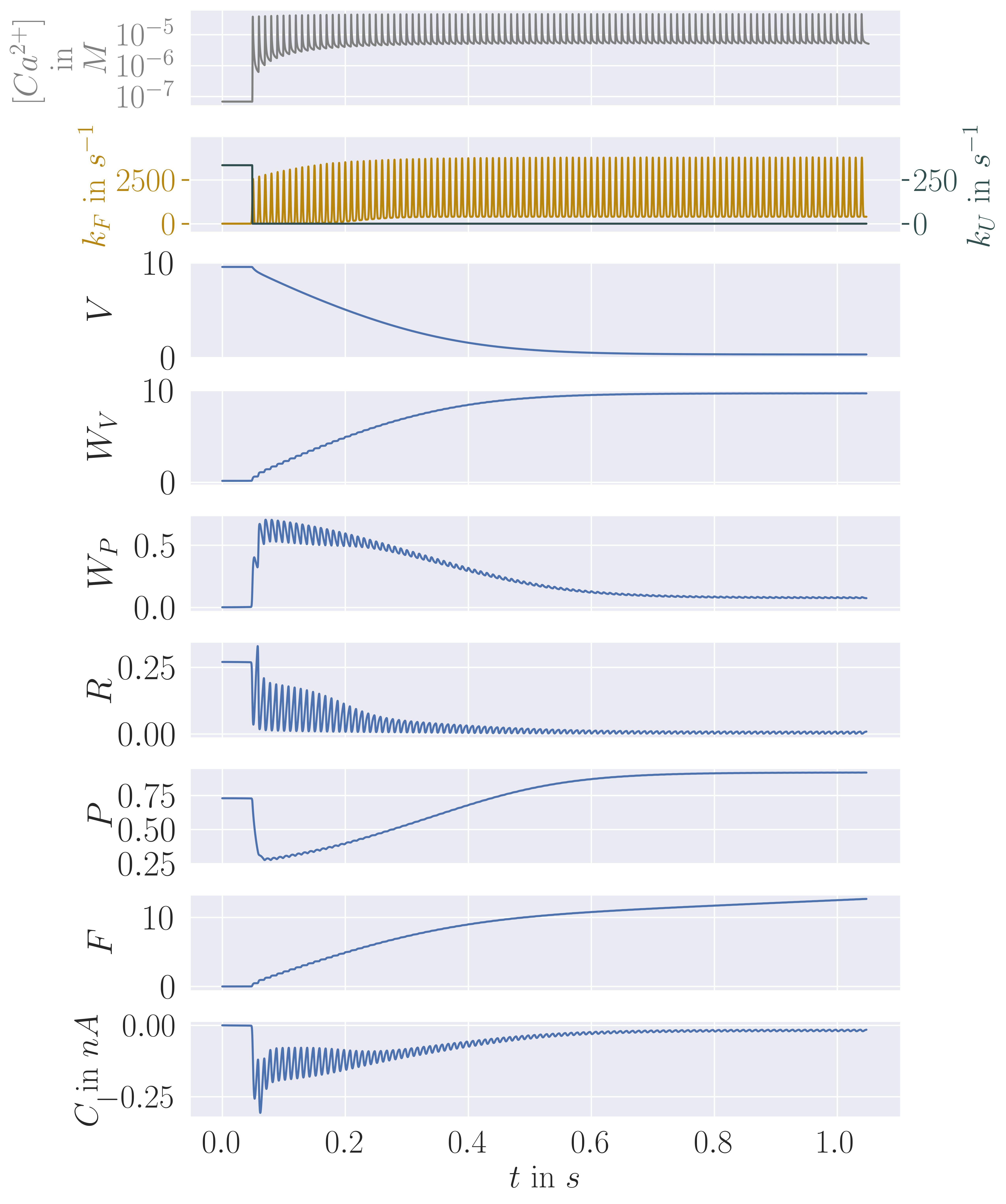}}
  \caption{\textbf{System response to sustained stimulation.} Temporal evolution of species' average numbers and the current $C$ (blue lines in third to ninth plot from the top, see respective labels on $y$-axes) in response to changes in the time-dependent rates ($k_F$ and $k_U$, shown as yellow and green lines in the second plot) representing a stimulus train of length $\SI{1}{\second}$ at frequency $f_{\text{stim}}=\SI{100}{\hertz}$, see Sec.~\ref{sec:sustained_stimulation}. Both rates depend on the intracellular $\text{Ca}^{2+}$ concentration dynamics (grey line in first plot), which were determined using the \textit{CalC} modeling tool \cite{matveev2002new} at a physiological external $\text{Ca}^{2+}$ concentration of $\SI{1.5}{\milli\Molar}$ and a distance of \SI{118}{\nm} from the $\text{Ca}^{2+}$ channel.  
 }
  \label{fig:normal}
\end{figure}

Fig.~\ref{fig:normal} demonstrates the temporal evolution of the model system and the current $C$ as response to sustained stimulation for $\SI{1}{\second}$ at frequency $f_{\text{stim}}=\SI{100}{\hertz}$, represented by the fusion rate $k_F(t)$ and the unpriming rate $k_U(t)$ (see Fig.~\ref{fig:normal}, second plot). Both functions depend on the intracellular $\text{Ca}^{2+}$ concentration dynamics (shown in the top plot of Fig.~\ref{fig:normal}), which were determined using the \textit{CalC} modeling tool \cite{matveev2002new} at a physiological external $\text{Ca}^{2+}$ concentration of $\SI{1.5}{\milli\Molar}$ and a distance of \SI{118}{\nm} from the $\text{Ca}^{2+}$ channel. Based on the $\text{Ca}^{2+}$ flow, we estimated the asymptotically periodic fusion rate $k_F$ as a weighted average of the fusion rates in the Kobbersmed model \cite{kobbersmed2020rapid} (see Sec.~\ref{parameter_estimation} for details). Also, following \cite{kobbersmed2020rapid}, we adapted the sigmoidal shape of $k_U(t)$, with the specific parameter values of this function and the priming rate constant $k_R$ chosen such that the facilitation effect was reproduced proportionally. For more details on the estimation of these rates and the remaining rate constants see Sec.~\ref{parameter_estimation}.

The numbers of release sites and vesicles were set to be $n_{\text{sites}}=1$ and $n_{\text{ves}}=10$, respectively, which means that we consider the average dynamics per release site assuming that the number of vesicles per release site is $10$. 
The system was initialized in steady state at $t=0$, i.e.,  $\boldsymbol{X}(0)=\hat{\boldsymbol{x}}$ is given by $h(\hat{\boldsymbol{x}},0)=0$, i.e., with reaction rates referring to no stimulation. The initial number of fusion events was set to $F(0)=0$. For the starting time of stimulation we chose $t_\text{start}=\SI{0.05}{\second} $.

The quantity analogous to experimentally measured currents is given by $C(t)$ shown in the bottom plot in Fig.~\ref{fig:normal}. The signal exhibits distinct phases: an initial large response (the first two stimuli, including a facilitation effect) and fast depression to a plateau lasting for about $\SI{0.1}{\second}$ and then a second slower decay to a final periodic orbit with significantly smaller amplitude. 
This behavior can be explained qualitatively: In the initial state (given by the steady state related to the initial parameter values), release sites are distributed between $R$ and $P$ due to the balance between the priming and unpriming reaction. The surplus of vesicles is accumulated in $V$. After the first stimulus, the unpriming rate drops to a very low value and release sites in $P$ can quickly bind to vesicles in $V$, which explains the initial facilitation and the strong response that is weakened quickly as release sites accumulate in $W_P$. Afterwards, the large vesicle supply in $V$ immediately provides recovered release sites with a vesicle and is thus gradually vacated while the signal plateaus. Due to the low vesicle recovery rate, vesicles start to accumulate in $W_V$. Once the amount of vesicles in $V$ approaches low values, increasingly more release sites are starting to collect in $P$ again and the system converges to a periodic orbit with a small signal amplitude.

\paragraph{Periodically forced system.} 

The final periodic orbit stems from sustained stimulation in which the rate $k_F$ depends on time periodically,
at least in an asymptotic sense, i.e., there is some time $t_0$ after which the dependence of $k_F$ can be considered periodic with period $T$ given by the stimulation frequency,
\[
k_F(t) = k_F(t+T),\quad \forall t\ge t_0,
\]
while $k_U$ is constant for $t\ge t_0$, and all other rates are time-independent. Consequently, the right-hand side function $h$ in \eqref{RRE} also is $T$-periodic via its dependence on $k_F$. In this case, dynamical systems theory \cite{Novaes2022,Capietto1992,Hausrath1988} tells us that, as long as the amplitude of the periodic forcing is not too large, there is at least one asymptotically $T$-periodic solution $\boldsymbol{X}_{\text{per}}=\boldsymbol{X}_{\text{per}}(t)$ of (\ref{RRE}) that oscillates around a certain fixed point $\boldsymbol{X}_0$. By continuation theory and averaging, we know that this fixed point $\boldsymbol{X}_0$ is the unique steady state of the averaged right-hand side function \cite{Novaes2022}
\begin{equation}\label{barf}
\bar{h}(\boldsymbol{X}):= \frac{1}{T}\int_{t_0}^{t_0+T} h(\boldsymbol{X},t) dt 
= \left(\begin{matrix}
-k_R V P +g_V W_V +\bar{k}_U R\\
\bar{k}_F R -g_V W_V\\
\bar{k}_F R - g_P W_P\\
k_R V P  -  \bar{k}_F R -\bar{k}_U R\\
- k_R V P +g_P W_P +\bar{k}_U R
\end{matrix}\right),
\end{equation}
with
\[
\bar{k}_F=\frac{1}{T}\int_{t_0}^{t_0+T} k_F(t) dt,\quad \bar{k}_U= k_U(t_0) .
\]
That is, $\boldsymbol{X}_0$ is the unique solution of $\bar{h}(\boldsymbol{X}_0)=0$, and, if $k_F(t)=\bar{k}_F+\lambda \tilde{k}_F(t)$ for  $T$-periodic $\tilde{k}_F$ with running average $0$, then we have that the periodic solution $\boldsymbol{X}_{\text{per}}$ converges to $\boldsymbol{X}_0$ for $\lambda\to 0$, and stays in a $\lambda$-wide neighborhood of $\boldsymbol{X}_0$ for not too large $\lambda$.

\begin{figure}[ht]
  \centering
   {\includegraphics[width=0.7\textwidth]{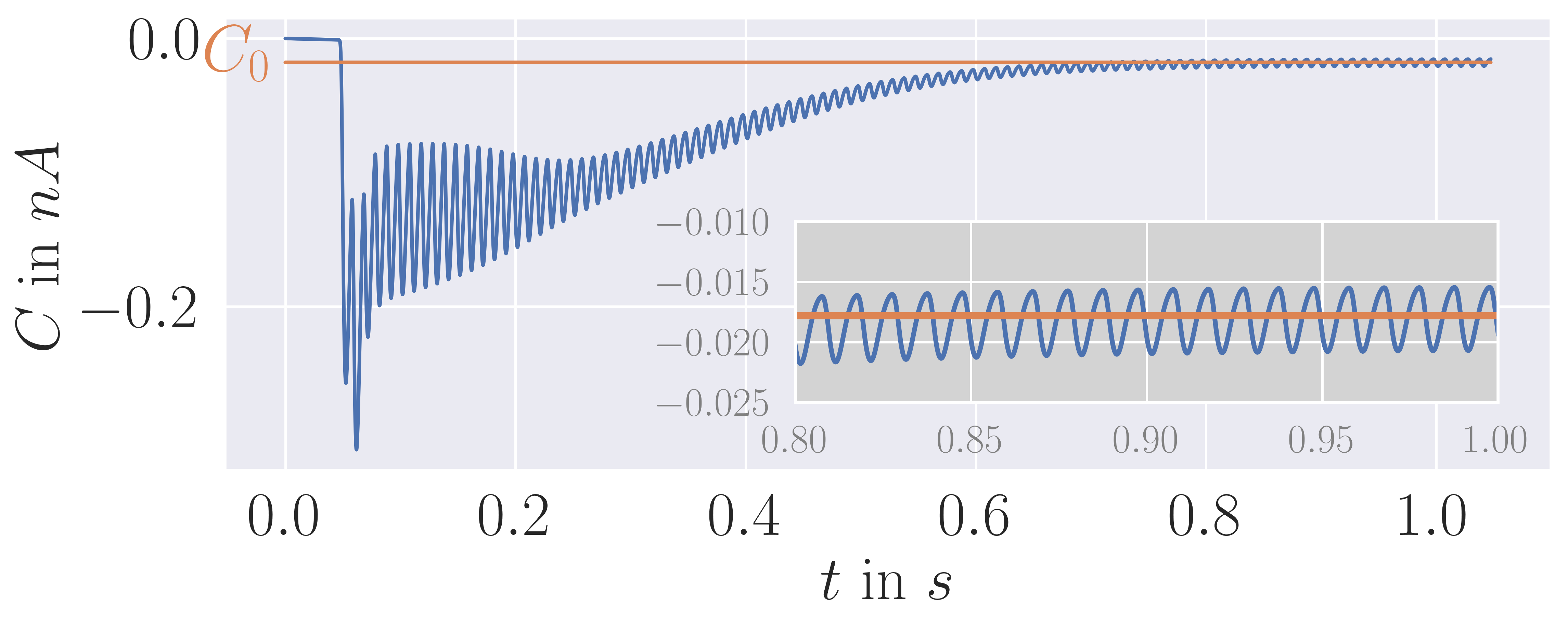}}
  \caption{\textbf{Asymptotic behavior.} Temporal evolution of the current $C$ (blue) and asymptotic oscillation around the fixed point $C_0$ (orange).  The grey inset shows a zoom-in of the last $\SI{0.2}{\second}$. The temporal averaging was done for $T=\frac{1}{\SI{100}{\hertz}}=\SI{0.01}{\second}$ and $t_0=\SI{0.99}{\second}$.
 }
  \label{fig:fixed_point}
\end{figure}

That is, for sustained periodic stimulation the system will asymptotically show periodic behavior with oscillations around the steady state given by the time-averaged rates. 
According to Eqs.~\eqref{eq:dFdt} and \eqref{output_signal}, this holds also true for the current $C$, which oscillates around the fixed point $C_0$ given by
\begin{equation}
   C_0 =  (\bar{k}_F R_0 * g)(t) = \bar{k}_F R_0 \int_{-\infty}^{\infty} g(s) ds.
\end{equation}
With the investigated parameter values, the asymptotic behavior can already be observed after less than $\SI{1}{\second}$ of stimulation, as depicted in Fig.~\ref{fig:fixed_point}.

\quad

The final periodic orbit as well as the transient system response are dictated by the balance between the two recovery processes. In the following we will consider the signal sensitivity in order to determine which process is more influential (i.e. 'rate-limiting' or 'rate-determining').

\subsection{Sensitivity analysis}

A characteristic property of the limiting process is that the signal $C(t)$ should be particularly sensitive to changes in its rate constant: for example, if vesicle recovery is more impactful than release site recovery, small changes in $g_V$ should result in a greater change in $C$ than small changes in $g_P$. We therefore introduce the notation $C(t,p(t))$ to emphasize that the output $C$ depends on the (partially) time-varying parameter values $p(t):=(k_R,k_U(t),k_F(t),g_V,g_P)$ and consider the \textit{sensitivity} $Z_{C}$  as a measure of the influence of the two recovery processes on the output signal $C$: 
\begin{equation}\label{eq:sens}
    Z_{C}^{g_V}(t) := \frac{\partial C}{\partial g_V}(t,p(t))\Big|_{p(t)=p^*(t)}, \quad \quad Z_{C}^{g_P}(t) := \frac{\partial C}{\partial g_P}(t,p(t))\Big|_{p(t)=p^*(t)},
\end{equation}
where $p^*$ refers to the parameter values given by the parameter estimation, see Sec.~\ref{parameter_estimation}. 
Defining the sensitivities $Z_{\dot F}^{g_V}(t)$, $Z_{\dot F}^{g_P}(t)$ in analogy to \eqref{eq:sens}, we observe that 
\begin{align}
  \frac{\partial }{\partial g_V}C(t)& = \frac{\partial}{\partial g_V} \int_{-\infty}^{\infty} \dot F(s)g(t-s) ds \\
  &\stackrel{(\star)}{=} \int_{-\infty}^{\infty}    \frac{\partial}{\partial g_V} \dot F(s) g(t-s) ds \\
  & = \int_{-\infty}^{\infty}  Z^{g_V}_{\dot F}(s) g(t-s) ds = (Z^{g_V}_{\dot F} \ast g)(t),
\end{align}
where we skipped $p(t)$ and used the Leibniz rule in $(\star)$. 
That is, we have 
\begin{equation}
    Z^{g_V}_C(t)=\left(Z^{g_V}_{\dot F} \ast g\right)(t),
\end{equation}
and analogously for $g_P$. Here, we used that the impulse response function $g$ does not depend on the parameter values $p$. 

The quantity $Z_{C}^{g_V}(t)$ captures the change in $C(t)$ induced by increasing the rate constant $g_V$ by an infinitesimal amount for all times, and likewise for $g_P$.
A closed system of ordinary differential equations can be derived for the sensitivities of all model components, which we can solve  simultaneously with the RRE in order to compute the sensitivities in $C$ at any time $t\geq 0$ \cite{dickinson_sensitivity_1976}. Further details can be found in Sec.~\ref{sec:sens_simple}. We finally normalize the sensitivity coefficients and define \cite{kirch2016effect}
\begin{equation}\label{sens_normalized}
    z^{g_V}_{C}(t):=Z^{g_V}_{C}(t)\cdot \frac{g_V}{C(t)}, \quad z^{g_P}_{C}(t):=Z^{g_P}_{C}(t)\cdot \frac{g_P}{C(t)}.
\end{equation}
Hereby, we obtain sensitivity values \textit{relative} to the rates $g_V,g_P$ and to the signal $C(t)$. This is especially important because $g_V$ is much smaller than $g_P$ in our parameter estimation, such that the absolute sensitivities $Z_{C}^{g_V}(t),Z_{C}^{g_P}(t)$ would deliver a distorted impression of the parameters' impact. 

\begin{figure}[ht]
  \centering
   {\includegraphics[width=0.8\textwidth]{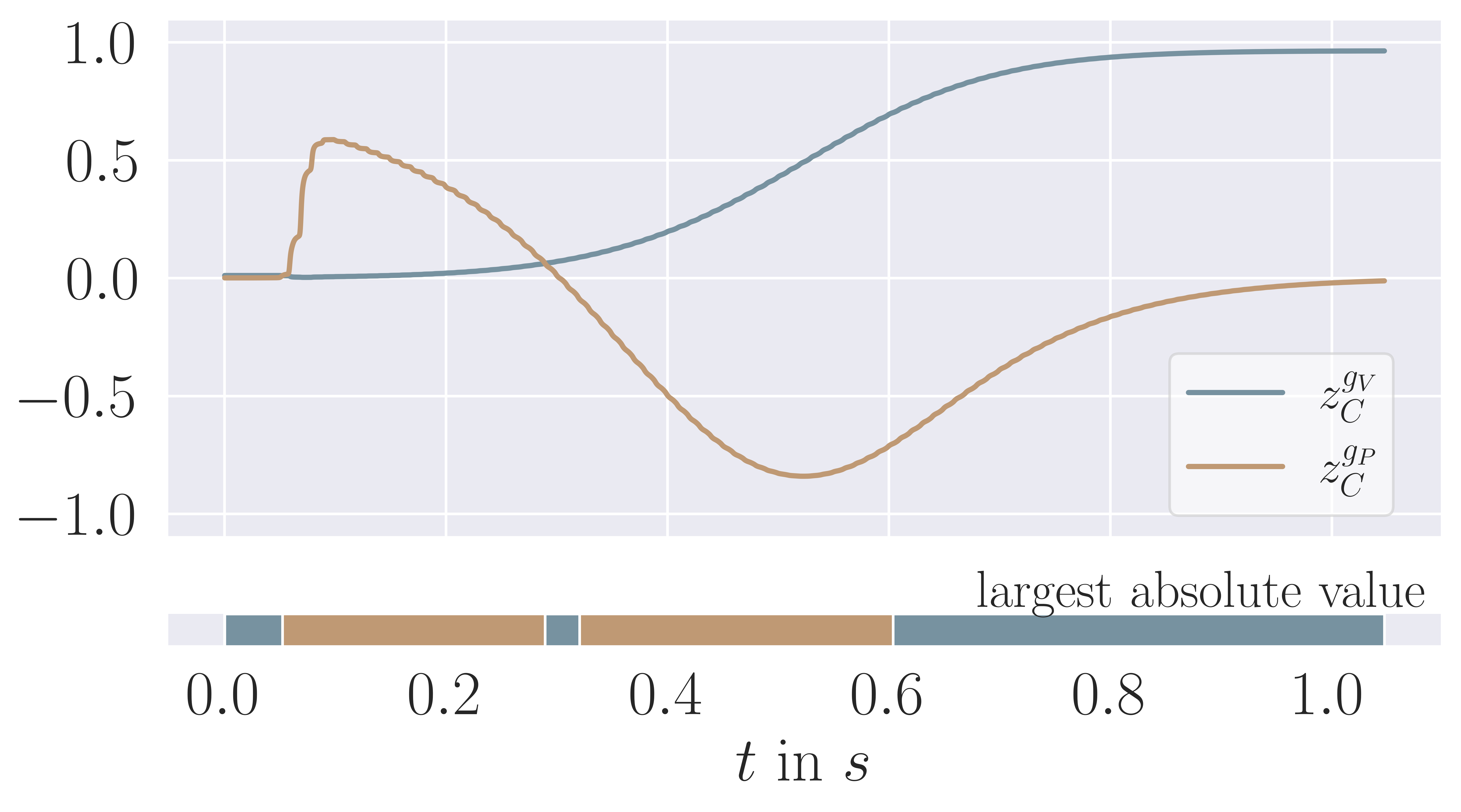}}
  \caption{\textbf{Normalized sensitivities under sustained stimulation.} Temporal evolution of the normalized sensitivity (given in Eq.~\eqref{sens_normalized}) of the current $C$ to vesicle and release site recovery rate. The system parameters were chosen as in Fig.~\ref{fig:normal}. The colored bar (bottom) indicates   which of the two sensitivities dominates in absolute values.}
  \label{fig:sens_normal}
\end{figure}

The temporal evolution of these normalized sensitivities for the system response pictured in Fig.~\ref{fig:normal} is shown in Fig.~\ref{fig:sens_normal}. 
At the very beginning, the sensitivities are almost zero for both recovery processes. This is because the system is initially at steady state without stimulation, where both $W_P$ and $W_V$ have very low numbers and recovery is of low significance for the resulting signal. The further temporal evolution in Fig.~\ref{fig:sens_normal} matches well with the above qualitative discussion on Fig.~\ref{fig:normal}: For $t\leq\SI{0.2}{\second}$ (during the plateau), the sensitivity to changes in $g_V$ (blue) is low and only slowly increasing due to the fact that most vesicles are stored in $V$. Once more vesicles start to accumulate in $W_V$, $z_{C}^{g_V}$ grows and finally approaches a constant positive value once the system reaches the asymptotic periodic orbit. 
The sensitivity in $g_P$ (orange) initially quickly rises 
because release site abundance increases in $W_P$  while it simultaneously decreases in $P$. As there are sufficient vesicles available in $V$ during the plateau for $t\leq \SI{0.2}{\second}$, faster release site recovery increases the signal and the sensitivity is positive. However, this also means that the supply in $V$ is emptied faster and the signal decreases at an earlier point in time. 
This is why the sensitivity $z_{C}^{g_P}$ becomes negative after the plateau at around $t=\SI{0.25}{\second}$. With the system settling into its final periodic orbit, $z_{C}^{g_P}$ approaches a constant low value. This is due to the fact that an increase in $g_P$ leads to a time-shift of the transient phase (from plateau to asymptotic orbit) to an earlier time period, but not to a significant change in the final orbit itself. 

In order to determine which of the recovery processes is more influential, one needs to compare the two sensitivities' absolute values (colored bar at the bottom of Fig.~\ref{fig:sens_normal}). Initially, for $t_\text{start}\leq t\leq \SI{0.2}{\second}$, the sensitivity $z_{C}^{g_P}$ with respect to $g_P$ clearly exceeds the value of $z_{C}^{g_V}$, which means that release site recovery is the limiting process. Near the shift from positive to negative values in $z_{C}^{g_P}$, the sensitivity with respect to $g_V$ temporally dominates, but then it is again surpassed by the negative impact that a permanent increase in the recovery rate $g_P$ has onto the signal in the time frame $\SI{0.3}{\second}\leq t \leq \SI{0.6}{\second}$. 
After time $t\approx \SI{0.6}{\second}$, with the system approaching the final periodic orbit, the sensitivity with respect to $g_V$ clearly dominates, while the impact of the recovery rate $g_P$ may be neglected. That is, in the long run, vesicle recovery is the limiting process.


This behavior leads us to an important general insight independent of the specific model and parameters: the answer to the question  of the rate-limiting process is not necessarily binary but can depend on the point in time during stimulation.


\paragraph{Dependence on parameter values.} Of course, the sensitivity evolution in Fig.~\ref{fig:sens_normal} is only the result for the specific set of parameter values estimated in Sec.~\ref{parameter_estimation}. However, parameter studies in which we varied either $k_R$, $g_V$, or $g_P$ from one twentieth to twenty times the original value demonstrate that the identity of the rate-limiting process is indeed a time-varying quantity over a large section of parameter space (details and images included in Sec.~\ref{sec:ParamStudies}). In all cases examined, upon stimulation, the dominant limiting process is initially release site recovery. The determining process then switches to vesicle recovery for some time unless $g_V$ is very large or $g_P$ is very small, that is, unless vesicle recovery is very fast by comparison. Afterwards, site recovery  may become limiting again for a period of time, but the system eventually switches back to higher sensitivity in $g_V$ within the $\SI{1}{\second}$ of stimulation in almost all cases (again, unless $g_P$ is very small). In summary, for a significant part of parameter space, the identity of the limiting process initially starts out as release site recovery but changes to vesicle recovery by the end of $\SI{1}{\second}$ of stimulation time, with the possibility of an additional switch to site recovery and back in between. This behavior can be ascribed to the initial surplus of vesicles in $V$ which leads to very low sensitivity to $g_V$. If $V$ is emptied and vesicles are accumulating in the recovery state, changes in $g_V$ have much higher impact and $z^{g_V}_C$ dominates. Again, this will happen unless vesicle recovery is very fast by comparison.


\quad 

By choosing $n_{\text{sites}}=1,n_{\text{ves}}=10$, the considered ODE is used to describe the \textit{average} dynamics at a \textit{single} release site which is available to $10$ vesicles. Augmenting the values $n_{\text{sites}}$ and $n_{\text{ves}}$ would mean to consider an active zone of several release sites which all access the same vesicle pool of size $n_{\text{ves}}$. Typically, there are only very few release sites per active zone which justifies to stick to small numbers $n_{\text{sites}}$, as done in this work. 
In this case, stochastic effects in the dynamics may play an important role, which motivates to extend the analysis to stochastic dynamics. 

\section{Stochastic dynamics of individual release sites}\label{sec:stochastic}

It is well-known that an ODE-description in terms of the reaction rate equation (as given by Eq.~\eqref{RRE}) delivers a good approximation of the average reaction dynamics in case of large particle numbers. For a single release site, however, the number of partaking vesicles is rather small and deviations from the ODE-behavior are to be expected. Furthermore, experimentally measured postsynaptic currents exhibit noise and irregularities even though multiple release sites are involved and summed over in the neurotransmission process. This motivates to investigate stochastic effects and variances of the recovery dynamics introduced in Sec.~\ref{sec:simplemodel} by formulating and analyzing the corresponding stochastic reaction-jump process.

\subsection{The reaction jump process} \label{sec:RJP}

The Markov process describing the stochastic recovery dynamics is denoted by 
\begin{equation}\label{X_cal}
    \boldsymbol{\mathcal{X}}(t)= (\mathcal{X}_i(t))_{i=1,...,5} = \left(\mathcal{V}(t),\mathcal{W}_V(t),\mathcal{W}_P(t),\mathcal{R}(t),\mathcal{P}(t)\right)^\top
\end{equation}
with $\boldsymbol{\mathcal{X}}(t) \in \mathbb{N}_0^5$ for all $t\geq 0$. Here, $\mathcal{X}_i(t)\in \mathbb{N}_0$ is the (random) number of vesicles or release sites in the respective states at time $t$, see again Fig.~\ref{fig:simplified_model} for the underlying model. These numbers change by discrete jumps induced by individual reaction events (given by the reactions \eqref{reaction1}-\eqref{reaction4}), which  occur after exponentially distributed waiting times. The associated probability distribution is characterized by the corresponding chemical master equation, see \cite{winkelmann2020stochastic} and references therein.

In analogy to Eq.~\eqref{output_signal}, the \textit{stochastic output current} is given by
\begin{equation}\label{C_cal}
        \mathcal{C}(t) := ( f \ast g)(t), 
\end{equation}
where  $g$ is again the impulse response function defined in Eq.~\eqref{eq:g} and $f$ is the functional derivative of the trajectories of the stochastic process $(\mathcal{F}(t))_{t\geq 0}$ counting the number of fusion events given by reaction \eqref{reaction3}. The latter is a monotonically increasing Markov jump process on the natural numbers, starting with $\mathcal{F}(0)=0$ and augmenting by one whenever a fusion event happens. Denoting the random jump times of $(\mathcal{F}(t))_{t\geq 0}$ by $T_1,T_2,...$, the functional derivative $f$ is a sum of Dirac delta functions shifted by the times $T_i$. 

\begin{figure}[ht]
		\centering
        \includegraphics[width=1.\textwidth]{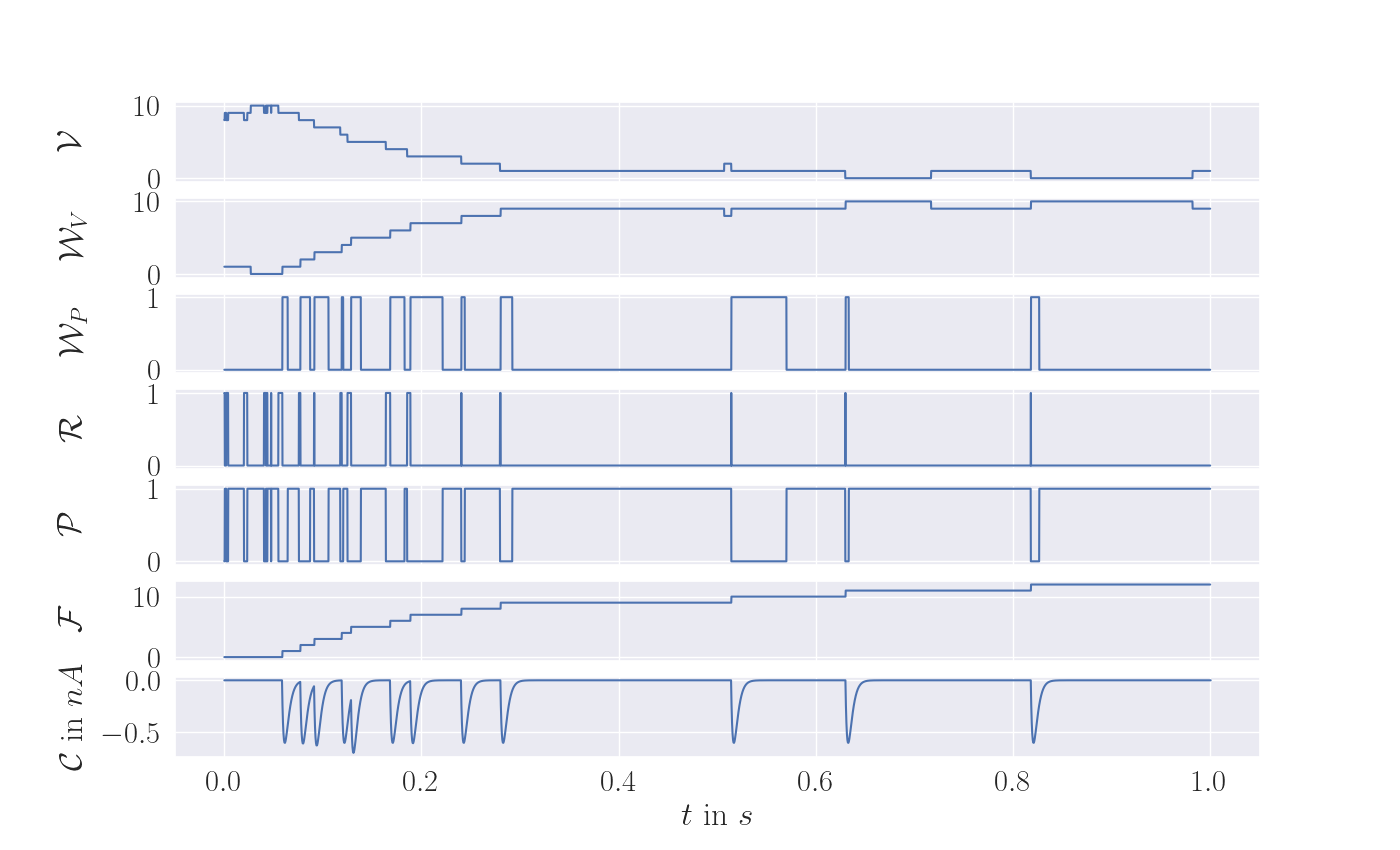}
		\caption{\textbf{Random dynamics at a single release site.} One realization of the stochastic reaction jump process $\mathcal{X}(t)$, the counting process $\mathcal{F}(t)$ and the stochastic output current $\mathcal{C}(t)$ for the same parameter values as used in Fig.~\ref{fig:normal}. The initial state is drawn randomly from the initial distribution which is the steady state distribution under no stimulation. 
  }
		\label{fig:Neuro_traj}
	\end{figure}

An individual trajectory of the stochastic reaction jump process $\mathcal{X}(t)$ is plotted in Fig.~\ref{fig:Neuro_traj}, including the counting process $\mathcal{F}(t)$ of fusion events and the induced stochastic output current $\mathcal{C}(t)$. 
This trajectory refers to the random dynamics at \textit{one single} release site. Its characteristics drastically deviate from the transient oscillatory and asymptotically periodic dynamics of the ODE-solution shown in Fig.~\ref{fig:normal}, although all parameter values coincide. The components $\mathcal{W}_P$, $\mathcal{R}$ and $\mathcal{P}$ switch between the discrete states $0$ and $1$, and the output current $\mathcal{C}$ consists of a few peaks occurring at random points in time. Especially, there is no obvious periodicity in the stochastic dynamics for such a single release site. However, the periodicity reappears by either considering the total junction current triggered by several release sites, which will be done in the following Sec.~\ref{sec:C_total}, or by calculating the dynamics' first-order moments, see Sec.~\ref{sec:moments}.

\subsection{Total junction current} \label{sec:C_total}
\begin{figure}[ht]
		\centering
		{\includegraphics[width=.9\textwidth]{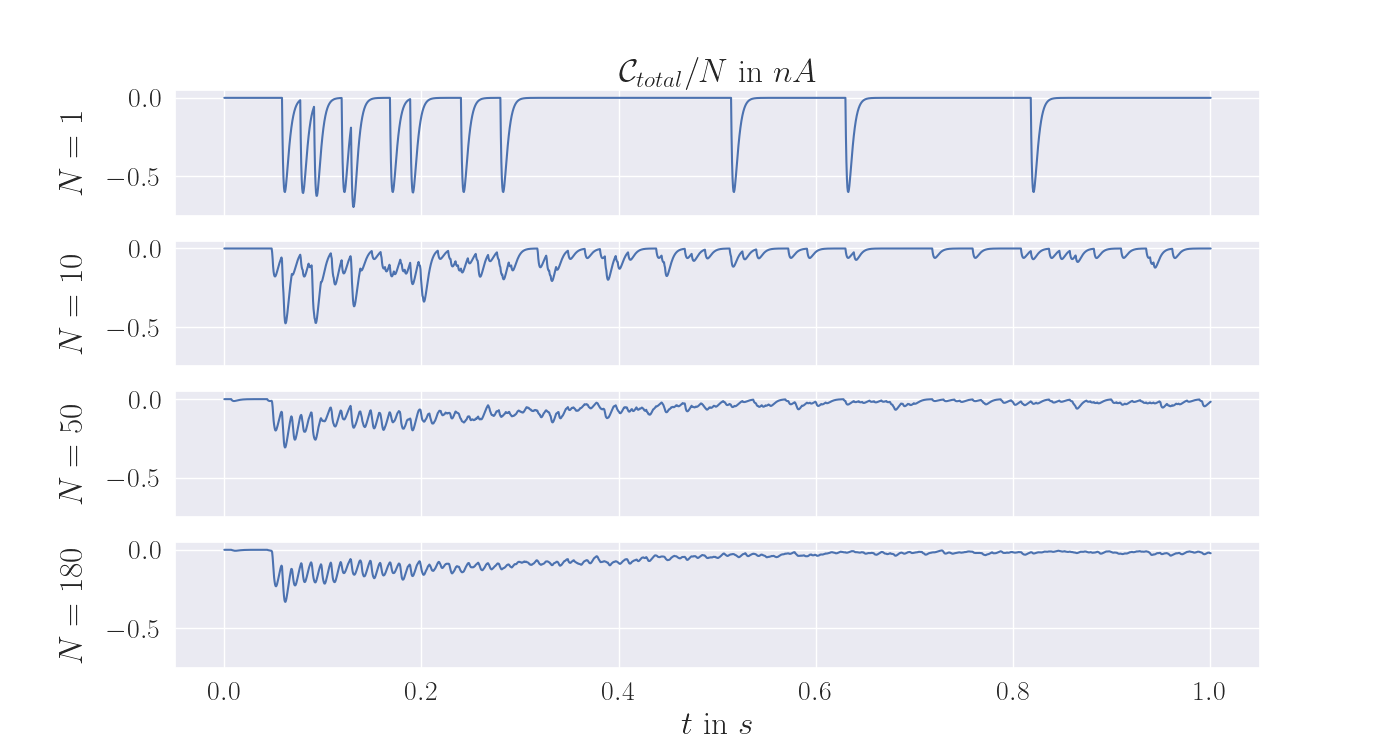}}
		\caption{\textbf{Scaled total output for different numbers of release sites.} One realization of the scaled total junction current $\mathcal{C}_{\text{total}}(t)/N$ for $N=1,10,50,180$.}  
		\label{fig:C_total_varyN}
	\end{figure}

Experimental measurements are typically given by the joint output signal of several release sites (in Ref.~\cite{ernst2022variance} we took $N=180$ release sites). The analogue quantity in our setting is given by the sum of $N\in \mathbb{N}$ independent 
realizations $\mathcal{C}_i$ of $\mathcal{C}$:
\begin{equation}\label{C_total}
    \mathcal{C}_{\text{total}}(t) = \sum_{i=1}^N \mathcal{C}_i(t).
\end{equation}
Fig.~\ref{fig:C_total_varyN} shows random realizations of the scaled total output current $\mathcal{C}_{\text{total}}/N$ for different numbers $N$ of release sites. 
For all $N$ one can observe that releases become scarcer in the course of time which is due to the fact that the reserve of vesicles is depleted. 
A periodicity in the dynamics is only perceptible for large $N$ ($N=50$, $N=180$) during the first $\SI{0.4}{\second}$ of stimulation. The characteristics pass from apparently non-periodic, randomly occurring peaks for small $N$ to periodic dynamics that appear to be close to the ODE-solution for large $N$. This can be explained as follows: By the law of large numbers, the scaled total output $\mathcal{C}_{\text{total}}(t)/N$ converges to the mean $\mu_{\mathcal{C}}(t)$ of $\mathcal{C}(t)$ which in turn is close to the ODE-solution, as we will show in the following Sec.~\ref{sec:moments}. For small $N$, the periodicity is hidden in the time-dependent fusion rate $k_F(t)$ and only becomes visible when looking at statistical averages. 
In general, the stochastic dynamics show large variations  which gradually decrease when increasing the number $N$ of release sites. The following investigation of the system's first- and second-order moments clarifies this issue.  


\subsection{First- and second-order moments}\label{sec:moments}

\begin{figure}[ht]
		\centering
  \begin{subfigure}[b]{1.\textwidth}
	    	\centering
    \includegraphics[width=1.\textwidth]{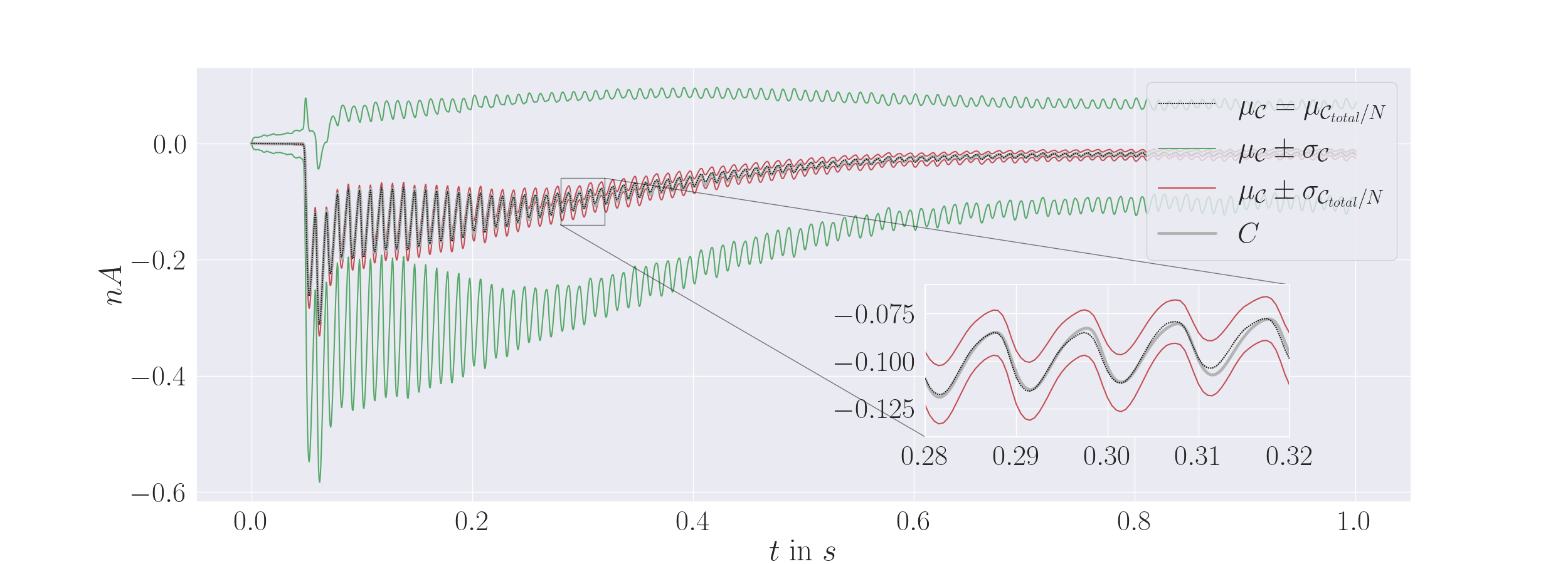}
	\end{subfigure}
		\caption{\textbf{First- and second order moments of the output current.} ODE-solution $C(t)$ (grey line) and first-order moment $\mu_{\mathcal{C}}(t)=\mu_{\mathcal{C}_{\text{total}}}(t)/N$ (black dotted line), together with the standard deviations $\sigma_{\mathcal{C}}(t)$ and $\sigma_{\mathcal{C}_{\text{total}}/N}(t)$ (green lines) of the stochastic output current $\mathcal{C}(t)$ and the rescaled total current $\mathcal{C}_{\text{total}}(t)/N$ for $N=180$, respectively. The zoom-in shows a high level of agreement between $C(t)$ (grey line) and $\mu_{\mathcal{C}}(t)$ (black dotted line). The relative standard deviation decreases substantially when considering the total current of several release sites. Estimated from $10^4$ MC simulations. Parameter values were chosen as in Fig.~\ref{fig:normal}. 
        }
		\label{fig:ave_std_zoom}
	\end{figure}

The first-order moment $\mu_{\mathcal{C}}(t):=\mathbb{E}(\mathcal{C}(t))=\frac{1}{N}\mathbb{E}(\mathcal{C}_{\text{total}}(t))$  of the stochastic signal $\mathcal{C}(t)$ and of the scaled total output current $\mathcal{C}_{\text{total}}(t)/N$ 
is plotted in Fig.~\ref{fig:ave_std_zoom}, together with the time-dependent standard deviations $\sigma_{\mathcal{C}}(t)$ and $\sigma_{\mathcal{C}_{\text{total}}/N}(t)$ for $N=180$, all estimated from $10^4$ MC simulations. We observe a periodicity in the first- and second-order moments as well as a close agreement of the mean $\mu_{\mathcal{C}}(t)$ with the ODE-solution $C(t)$ from Fig.~\ref{fig:normal}. This similarity is surprising because nonlinear reaction systems typically show a significant deviation of the stochastic mean from the ODE-solution, at least for small particle numbers \cite{winkelmann2020stochastic}, which would imply the inequality $\mu_{\mathcal{C}}(t) \neq C(t)$.
However, further numerical experiments on the reaction system under investigation 
show that the high-level similarity $\mu_{\mathcal{C}}(t) \approx C(t)$ exists independently of the population size and the chosen parameter values. Indeed, the source of the similarity mainly lies in the independence of the recovery processes of vesicles and release sites which implies a small covariance of their dynamics, as we will explain in the following.  


\begin{figure}[ht]
			\centering
	\includegraphics[width=1.\textwidth]     {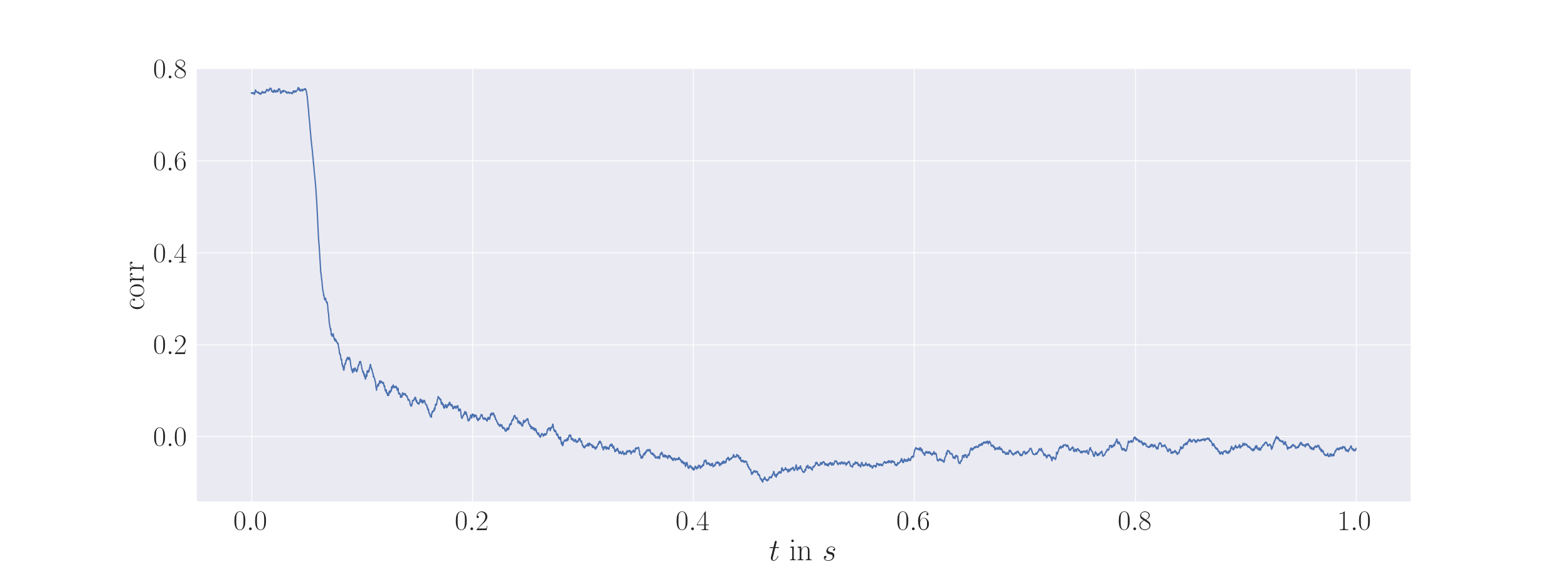}
	   \caption{\textbf{Temporal evolution of the correlation between $\mathcal{V}(t)$ and $\mathcal{P}(t)$.} The correlation function $\text{corr}(t)$ is defined in Eq.~\eqref{eq:cor} and was estimated from $10^4$ MC simulations. Parameter values were chosen as in Fig.~\ref{fig:normal}.}
        \label{fig:correlation}
	\end{figure}

\paragraph{Independent recovery processes.} After a fusion event (which carries both the vesicle and the release site into their recovery states $W_V$ and $W_P$, respectively), the recovery dynamics given by reaction~\eqref{reaction4} happen independently of each other. That is, the time it takes a recovering vesicle to increase the number of available vesicles $V$ again does not affect the time it takes a recovering release site to add to the number of available sites $P$, and vice versa. This stands in contrast to the unpriming reaction~\eqref{reaction2}, which simultaneously augments both $V$ and $P$. However, unpriming  happens at a very small rate (after a short initial phase), so that the increase in $V$ or $P$ mainly results from independent recovery reactions. Thereby, we obtain a certain degree of independence in the distributions of $V$ and $P$, meaning that we have (with respect to the law of the jump process),
\begin{equation}
    \mathbb{P}\left[\mathcal{V}(t)=n,\mathcal{P}(t)=m\right] \approx \mathbb{P}\left[\mathcal{V}(t)=n\right] \cdot \mathbb{P}\left[\mathcal{P}(t)=m\right],
\end{equation}
for $n,m\in \mathbb{N}_0$ and $t> 0$, 
as well as
\begin{equation}\label{approx_E}
    \mathbb{E}\left[\mathcal{V}(t)\cdot\mathcal{P}(t)\right] \approx \mathbb{E}\left[\mathcal{V}(t)\right] \cdot \mathbb{E}\left[\mathcal{P}(t)\right]
\end{equation}
for most times $t\geq 0$, and consequently $\mathbb{E}(\boldsymbol{\mathcal{X}}(t)) \approx \boldsymbol{X}(t)$ and $\mathbb{E}(\mathcal{C}(t)) \approx C(t)$ . 
This similarity is equivalent to a small covariance 
\begin{equation}
    \text{cov}(\mathcal{V}(t),\mathcal{P}(t)) :=  \mathbb{E}\left[\mathcal{V}(t)\cdot\mathcal{P}(t)\right]- \mathbb{E}\left[\mathcal{V}(t)\right] \cdot \mathbb{E}\left[\mathcal{P}(t)\right]
    \label{eq:cov}
\end{equation} 
or a small correlation
\begin{equation}
    \text{corr}(t) := \frac{\text{cov}(\mathcal{V}(t),\mathcal{P}(t))}{\sigma_\mathcal{V}(t)\sigma_\mathcal{P}(t)},
    \label{eq:cor}
\end{equation}
where $\sigma_\mathcal{V}(t)$ and $\sigma_\mathcal{P}(t)$ denote the standard deviations of the processes $\mathcal{V}(t)$ and $\mathcal{P}(t)$, respectively. The correlation function is plotted in  Fig.~\ref{fig:correlation}. One observes a rapid decrease of the correlation towards zero, which corresponds to an increase in the degree of independence in our system.
We study the correlation for two reduced reaction systems in order to clarify this effect of the independent recovery dynamics in Sec.~\ref{reduced_reaction_system}. 
	

\section{Concluding Remarks}

In this work, we have introduced a non-linear reaction network for the presynaptic dynamics of signal processing 
at chemical synapses, including explicit reactions for recovery processes. 
Modeled by a second-order reaction,  a freely available vesicle attaches to a freely available release site. Afterwards, the vesicle either detaches again or fuses with the membrane, thereby triggering an output current. Both the unpriming rate and the fusion rate depend on time, accounting for the level of stimulation. After a fusion event, both the vesicle and the release site enter recovery states, where they stay until they return independently of each other and become available again. The goal of this work was (1) to understand the effect of the recovery processes on the total signal output under sustained stimulation of the system, and (2) to investigate how single release site dynamics may be related to the output current measured for many release sites as typically in experiments.  

We have analyzed the signaling process by numerically solving the associated reaction rate equation and by simulating the associated stochastic reaction jump process. The main findings may be summarized as follows. 
\begin{itemize}
    \item During the initial phase of stimulation, there is a relatively stable response caused by the combined effect of fast release site recycling and a high supply of vesicles available for binding to the sites. With sustained activity, the vesicle supply is depleted which leads to a deep depression of the output signal. This is caused by the small recovery rate of the vesicles. Finally, the signal reaches a periodic orbit with a small amplitude in the vicinity of a uniquely determined steady state of the averaged system.
    \item Sensitivity analysis for the recovery rates reveals a considerable time-dependence of the normalized sensitivity coefficients. While the output current is dominantly influenced by the recovery rate $g_P$ of release sites at the start of stimulation, at later times, the vesicle recovery rate $g_V$ becomes more decisive. This result is in good agreement with previous discussions on the subject \cite{neher2010rate} and extends them by the idea that the identity of the rate-determining step of neurotransmission under sustained activity depends on the time during stimulation. Since the evolution of the sensitivities depends on the parameter values, we conducted a parameter study and observed that the characteristic structure is conserved over a large range of values.  
    \item By simulating single release sites by means of  the stochastic reaction jump process we uncovered that the characteristics of an individual random trajectory, which contains separate spikes of the output signal occurring at random points in time, drastically deviate from the oscillatory and asymptotically periodic ODE-trajectory. At the same time, the first-order moments of the random dynamics showed a surprisingly close agreement with the ODE-solution. We have traced this closeness back to small correlation between vesicle supply and release site supply, which results from the near independence of the recovery processes. 
    Moreover, we found that the periodic pattern of the ODE-solution is recovered for the stochastic dynamics when one considers the total junction current induced by averaging over larger numbers of release sites. 
\end{itemize}


Overall, the model introduced in this paper allows to study the effects of recovery processes within presynaptic neurotransmission dynamics on different levels. The results presented herein guide the way towards qualitatively and quantitatively understanding the role of different recovery rates and stimulation frequencies on the overall output signal.
In this respect central questions of interest for future investigations include the following: Is it possible to differentiate between the effects caused by release site recovery in contrast to vesicle recovery? More concretely, may the same postsynaptic current arise from different combinations of recovery speeds? This question is of interest because it reveals to which extent experimental data may be used to uniquely identify the recovery rates. Answering this question would mean to quantitatively compare the system's response at different parameter values and stimulation frequencies and compare them with experimental measurements. 

\subsubsection*{Acknowledgements}
This research has
been partially funded by the Deutsche Forschungsgemeinschaft (DFG, German Research
Foundation) through grant CRC 1114/3 and under Germany’s Excellence Strategy – The
Berlin Mathematics Research Center MATH+ (EXC-2046/1 project ID: 390685689).

\subsubsection*{Code availability}
The code used to generate the results in this paper is available at \url{https://doi.org/10.5281/zenodo.7551439}.

\appendix

\section{Appendix}

\subsection{Sensitivity equation for recovery model}\label{sec:sens_simple}

Given the recovery model of Sec.~\ref{sec:simplemodel} we consider the extension
\begin{equation}
    \boldsymbol{Y}(t):=\left(V(t),W_V(t),W_P(t),R(t),P(t),F(t)\right)^\top
\end{equation}
of the cumulative state $\boldsymbol{X}(t)$ defined in \eqref{X}, which satisfies the extended RRE of the form
\begin{equation}\label{dotY}
   \Dot{\boldsymbol{Y}}(t) = 
\left(\begin{matrix}
-k_R V(t) P(t) +g_V W_V(t) \\
k_F(t) R(t) -g_V W_V(t)\\
k_F(t) R(t) - g_P W_P(t)\\
k_R V(t) P(t)  -  k_F(t) R(t)\\
- k_R V(t) P(t) +g_P W_P(t)\\
k_F(t) R(t)
\end{matrix}\right) =:\tilde h(\boldsymbol{Y}(t),t).
\end{equation}
We emphasize the dependency on the time-varying parameters $p(t)=(k_R,k_F(t),g_V,g_P)$ by writing $\boldsymbol{Y}(t,p(t))$ and $\dot{\boldsymbol{Y}}(t,p(t))$. 

We introduce the sensitivities 
\begin{equation}\label{eq:sensitivities_introduction}
    Z_{Y_i}^{g_V}(t):=\frac{\partial Y_i}{\partial  g_V}(t,p(t))\Big|_{p(t)=p^*(t)}, \quad   Z_{Y_i}^{g_P}(t):=\frac{\partial Y_i}{\partial  g_P}(t,p(t))\Big|_{p(t)=p^*(t)},
\end{equation}
with $p^*(t)$ being the reference parameter set from the parameter estimation, see Sec.~\ref{parameter_estimation}.
Let
\begin{align*}
\boldsymbol{Z}^{g_V}(t) &:= 
\left(\begin{matrix}
Z_V^{g_V}(t) \\
Z_{W_V}^{g_V}(t)\\
Z_{W_P}^{g_V}(t)\\
Z_{R}^{g_V}(t)\\
Z_{P}^{g_V}(t)\\
Z_{F}^{g_V}(t)
\end{matrix}\right), 
\hspace{.5cm}
\boldsymbol{Z}^{g_P}(t) := 
\left(\begin{matrix}
Z_V^{g_P}(t) \\
Z_{W_V}^{g_P}(t)\\
Z_{W_P}^{g_P}(t)\\
Z_{R}^{g_P}(t)\\
Z_{P}^{g_P}(t)\\
Z_{F}^{g_P}(t)
\end{matrix}\right).
\end{align*}
In accordance with \cite{dickinson_sensitivity_1976}, the ODE-systems for the sensitivities in matrix notation are then
\begin{align}
    \Dot{\boldsymbol{Z}}^{g_V}(t)  &= \frac{\partial \Dot{\boldsymbol{Y}} }{\partial g_V}(t,p(t))\Big|_{p(t)=p^*(t)} +\mathcal{J}\big(t,p^*(t)\big) \boldsymbol{Z}^{g_V}(t), \\
        \Dot{\boldsymbol{Z}}^{g_P}(t)  &= \frac{\partial \Dot{\boldsymbol{Y}} }{\partial g_P}(t,p(t))\Big|_{p(t)=p^*(t)} +\mathcal{J}(t,p^*(t)) \boldsymbol{Z}^{g_P}(t),
\end{align}
where \begin{align}
    \frac{\partial  \Dot{\boldsymbol{Y}}}{\partial g_V}(t,p(t)) &= \left(\begin{matrix}
    W_V(t,p(t))\\
    -W_V(t,p(t))\\
    0\\
    0\\
    0\\
    0
    \end{matrix}\right), \hspace{1cm }
    \frac{\partial  \Dot{\boldsymbol{Y}}}{\partial g_P}(t,p(t)) = \left(\begin{matrix}
    0\\
    0\\
    -W_P(t,p(t))\\
    0\\
    W_P(t,p(t))\\
    0
    \end{matrix}\right),
\end{align}
while $\mathcal{J}(t,p(t))$ denotes the Jacobian matrix of $\Dot{\boldsymbol{Y}}(t,p(t))$ with respect to $\boldsymbol{Y}$, so $\mathcal{J}_{ij}(t,p(t)):=\frac{\partial \Dot{Y}_i}{\partial Y_j}(t,p(t))=\frac{\partial \tilde h}{\partial Y_j}(\boldsymbol{Y},t)$ for $\dot{\boldsymbol{Y}}$ given in \eqref{dotY},  such that 
\begin{align}
  \mathcal{J}(t,p(t)) &= 
    \left( \begin{matrix}
    -k_R P(t,p(t)) &  g_V & 0 & 0 & -k_R V(t,p(t)) & 0 \\
    0 & -g_V & 0 & k_F(t) & 0 & 0 \\
    0 & 0 & -g_P & k_F(t) & 0 & 0 \\
    k_R P(t,p(t)) & 0 & 0 & -k_F(t) & k_R V(t,p(t)) & 0 \\
    -k_R P(t,p(t)) & 0 & g_P  & 0 & -k_R V(t,p(t))& 0\\
    0 & 0 & 0 & k_F(t)  & 0 & 0
    \end{matrix} \right).
\end{align}

We can solve numerically for the sensitivities of all model components. As the system is assumed to start in the parameter-dependent steady state $\hat{\boldsymbol{x}}=\hat{\boldsymbol{x}}(p)$ (see \ref{app:steady_state} for its calculation), the initial values of the sensitivities are given by 
\begin{equation}
Z^{g_P}_{X_i}(0)=\frac{\partial \hat{x}_i}{\partial g_P}(p(0))\Big|_{p(0)=p^*(0)}, \quad Z^{g_V}_{X_i}(0)=\frac{\partial \hat{x}_i}{\partial g_V}(p(0))\Big|_{p(0)=p^*(0)} \quad \mbox{for} \quad  i = 1,...,5
\end{equation}
for $(X_i)_{i=1,...,5}=(V,W_V,W_P,R,P)$. Moreover, as $F(0)=0$ holds independently of the parameter values, we know that $Z_{F}^{g_P}(0)=Z_{F}^{g_V}(0)=0$.

Finally, the sensitivity in $\Dot{F}$ can then be found by interchanging the partial derivatives, where one has to apply  Schwarz's theorem:
\begin{align*}
    Z^{g_V}_{\Dot{F}}(t) &= \frac{\partial \Dot{F}}{\partial g_V}(t) =  \frac{\partial}{\partial g_V} \frac{\partial F}{\partial t}(t) =  \frac{\partial }{\partial t} \frac{\partial F}{\partial g_V}(t) =  \frac{\partial }{\partial t} Z^{g_V}_{F}(t).
\end{align*}
The sensitivity $Z^{g_P}_{\Dot{F}}(t)$ can be found in an analogous manner.

\subsection{Estimation of  parameter values}\label{parameter_estimation}

In order to yield realistic model behavior, the rate constants were estimated both from the literature and from the \textit{unpriming model} by Kobbersmed et al. \cite{kobbersmed2020rapid}. In this model, there are several states that each release site can attend:
It can be either empty (state $P_0$), or there is a vesicle attached to it, which itself has zero to five $\text{Ca}^{2+}$ ions bound to its fusion sensor (states $R_0,...,R_5$, respectively). Switches between these states happen by random jump events, where the jump propensities partially depend on the $\text{Ca}^{2+}$ concentration. 
Besides the \textit{priming} reaction $P_0\to R_0$, which represents the binding of a vesicle to the empty release site, there is the reverse \textit{unpriming} reaction, which describes the process of a docked vesicle detaching from the release site again. 
The central event of a docked vesicle fusing with the membrane can happen from each of the states $R_0,...,R_5$ and turns the release site into the empty state $P_0$ again. This is modeled by the reaction $R_m\to P_0+F$ ($m=0,...,5$), where $F$ refers to the cumulative number of fusion events.

Our model, as described in Sec.~\ref{sec:simplemodel}, merges the states $R_0,...,R_5$ of Ref.~\cite{kobbersmed2020rapid} to one state $R$. On the other hand, we add the recovery states $W_V$ and $W_P$ as well as the state $V$ of available vesicles, thereby turning the first-order priming reaction $P_0\to R_0$ of the Kobbersmed model into a second-order reaction $P+V\to R$ between available release sites and available vesicles.  
These relations between the models are used in the following to estimate the fusion rate as well as the priming and unpriming rates for our model. 
An overview of the parameter values as well as the used method of estimation is shown in Table \ref{tab:parameters}, the details will be discussed in the following.

  \quad

 \noindent\textbf{Release site recovery rate $g_P\gtrsim \SI{50}{\per \second}$.}
    According to Kawasaki et al. \cite{kawasaki_fast_2000}, repeated stimulation in mutants with inhibited vesicle recovery induced synaptic fatiguing within $\SI{20}{\milli\second}$. Thus, release site recovery is estimated to  operate at a rate of $g_P=\frac{1}{\SI{20}{\milli\second}}=\SI{50}{\per \second}$. 
    
    \quad

    \noindent \textbf{Vesicle recovery rate $g_V\gtrsim \SI{0.4}{\per \second}$.} Watanabe et al. \cite{watanabe_clathrin_2014} observed and timed a succession of steps for vesicle recovery: endocytosis of a large vesicle ($\SI{50}{\milli\second}-\SI{100}{\milli\second}$), transition to an endosome ($\SI{1}{\second}$), coating ($\SI{3}{\second}$) and separation of the endosome into approximately 4 synaptic vesicles ($\SI{6}{\second}$). We therefore estimated the vesicle recovery rate to be 
    $g_V=[(\SI{100}{\milli\second} + \SI{1}{\second} + \SI{3}{\second} + \SI{6}{\second})/4]^{-1} \approx (\SI{2.5}{\second})^{-1} = \SI{0.4}{\per\second}$.  

    \quad
    
    
   \noindent 
    \textbf{Fusion rate $k_F(t)$.}  We estimated the fusion reaction propensity $k_F(t)$ during stimulation by calculating a weighted average of the fusion rates in the Kobbersmed model \cite{kobbersmed2020rapid} for $\SI{1}{\second}$  of stimulation at $\SI{100}{\hertz}$  with an external $\text{Ca}^{2+}$ concentration of $\SI{1.5}{\milli\Molar}$ and a distance from the $\text{Ca}^{2+}$ channel of \SI{118}{\nm}. 
    The weights result from truncating the states $P_0$ and $F$ from the model and observing the distribution of the release sites $R_0,...,R_5$ in the truncated model in response to the stimulus train.  In the Kobbersmed model, which is based on the allosteric fusion model by Lou et al. \cite{lou2005allosteric}, the dynamic behavior results from time-dependent changes in the intracellular $\text{Ca}^{2+}$ concentration that directly enters the model's reaction rates. The behavior of this $\text{Ca}^{2+}$ concentration was determined using the $CalC$ modeling tool \cite{matveev2002new} in accordance with the stimulation frequency and the number of applied stimuli (see top plot in Fig.~\ref{fig:normal}).

   In order to let $k_F$ be a continuously differentiable function, we approximate the resulting time-dependent weighted average $\bar{k}_F^{\text{Kob}}(t)$ as the sum of a baseline rate $f_{\text{baseline}}$ and a number of Gaussians $f_{\text{Gaussians}}$:
        \begin{align}
       k_F(t) &= f_{\text{baseline}}(t) + f_{\text{Gaussians}}(t),
    \end{align}
    where
    \begin{align}
     f_{\text{baseline}}(t)&=\frac{m_0}{1+e^{-m_1(t-m_2)}},\\
     f_{\text{Gaussians}}(t)&= \sum_{i=1}^{100}a_i e^{-\frac{0.5(t-t_{\text{stim},i})^2}{\sigma^2}},
    \end{align}
    with $m_0,m_1,m_2, \sigma, t_{\text{stim},i}, a_i \in \mathbb{R}_+$ for $i =1,..,100$, see Table~\ref{tab:parameters} for the values.
    The parameter values for the baseline function were found by fitting $f_{\text{baseline}}$ to the troughs of the weighted average $\bar{k}_F^{\text{Kob}}$. The parameter $m_0$ denotes the supremum of the logistic function $f_{\text{baseline}}$, while $m_1$ regulates the steepness and $m_2$ the time at which it assumes its midpoint.
    The peak times of the fusion rate $t_{\text{stim},i}$ and amplitudes $a_i$ with respect to the troughs were taken directly from $\bar{k}_F^{\text{Kob}}$, while the peak width $\sigma$ was approximated as the average peak width in $\bar{k}_F^{\text{Kob}}$. The resulting function $k_F(t)$ is plotted in the upper panel of Fig.~\ref{fig:normal}.

    \quad
    
    
    
   \noindent  \textbf{Priming rate $k_R$ and unpriming rate $k_U(t)$.} 
   In order to preserve the  paired-pulse ratio from the Kobbersmed model, we optimized both the time-dependent unpriming rate $k_U(t)$ and the priming rate $k_R$ in the follwing way. 
   In the interest of keeping the number of optimization parameters low, we assumed that $k_U(t)$ was of the same general shape as in Ref.~\cite{kobbersmed2020rapid} which can approximately be described with the following continuously differentiable sigmoid function:
        \begin{align}\label{kU}
        k_U(t) &= k^{\text{max}}_{U} \left(1-\frac{1}{1+e^{-m_3(t-m_4)}}\right) + k^{\text{min }}_{U}.
    \end{align}
       The parameters $m_3$, $m_4$, $k^{\text{min}}_{U}$ $\in \mathbb{R}_+$  were estimated directly by fitting this function to the unpriming rate from the Kobbersmed model and adopting the same values, see Table~\ref{tab:parameters}. The remaining parameters $k_R$, $k^{\text{max}}_{U} \in \mathbb{R}_+$ were then found by minimizing the parameter-dependent loss function
             \begin{equation}
        L(p) =  
         \abs*{\frac{\dot{F}(t_{\text{peak},2},p)}{\dot{F}(t_{\text{peak},1},p)} - \frac{\dot{F}^{\text{Kob}}(t_{\text{peak},2}^{\text{Kob}},p)}{\dot{F}^{\text{Kob}}(t_{\text{peak},1}^{\text{Kob}},p)}},
    \end{equation}
fixing the previously estimated parameter values of $p(t)=(k_R,k_U(t),k_F(t),g_V,g_P)$.
    Here, $t_{\text{peak},j}$ denotes the point in time of the $j$-th peak in $\dot{F}$ in the recovery model\footnote{Note that these times are not the same as $t_{\text{stim},i}$ from the previous section since there is some latency behind the rise of the fusion rate $k_F$ and the evocation of a signal.} and, accordingly, $t_{\text{peak},j}^{\text{Kob}}$ is the point in time of the $j$-th peak  in the unpriming model with $j=1,2$. The total number of vesicles and release sites were set to $n_{\text{ves}}=10$ and $n_{\text{sites}}=1$, respectively.
    
  \quad

  All parameter values are listed in Table~\ref{tab:parameters} and \ref{tab:ai} and determine the reference values $p^*(t)$ of the rate functions. 
    
    \begin{table}[H]\renewcommand{\arraystretch}{1.2}
\begin{tabular}{l|l|l}
\textbf{parameter} & \textbf{value} & \textbf{method of estimation} \\ \hline
$t_{\text{start}}$ & \SI{0.05}{\second} & chosen freely \\ \hline
$g_P$ & 50$s^{-1}$ & literature  \cite{kawasaki_fast_2000} \\ \hline
$g_V$ & 0.4$s^{-1}$ & literature \cite{watanabe_clathrin_2014} \\ \hline
$m_0$ &  \SI{397}{\per \second} & fitting $f_{\text{baseline}}$ to troughs of $\bar{k}_F^{\text{Kob}}$ \\
$m_1$ &  \SI{33.3}{\per \second}& ---''--- \\
$m_2$ & $t_{\text{start}} +$ \SI{0.174}{\second} & ---''--- \\
$a_i$ & see Table \ref{tab:ai}  & peak amplitudes with respect to troughs from $\bar{k}_F^{\text{Kob}}$ \\
$t_{\text{stim},i}$ & $t_{\text{start}}+ i\cdot$\SI{0.01}{\second}  & peak times from $\bar{k}_F^{\text{Kob}}$ \\
$\sigma$ & \SI{9.53e-4}{\per \second} & avg. peak width from $\bar{k}_F^{\text{Kob}}$ \\ \hline
$m_3$ & \SI{27318}{\per \second} & fitting $k_U$ to $k_U^{\text{Kob}}$ \\
$m_4$ & $t_{\text{start}} -$\SI{1.4e-3}{\second} &  ---''--- \\
$k^{\text{min}}_{U}$ & \SI{1.02e-08}{\second}  & ---''---\\
$k^{\text{max}}_{U}$ &   \SI{334}{\per \second}& minimizing $L(F(t))$ \\ \hline
$k_R$ & \SI{12.9}{\per \second} & minimizing $L(F(t))$
\end{tabular}\caption{\textbf{Parameter estimation results.} Overview of parameter values and the used method of estimation.}\label{tab:parameters}
\end{table}

\begin{table}[ht]
\begin{tabular}{l l||l l||l l||l l||l l}
$i$ & $a_i$ in $\SI{}{\per\second}$  & $i$  & $a_i$ in $\SI{}{\per\second}$ &$i$  & $a_i$ in $\SI{}{\per\second}$  &$i$  &$a_i$ in $\SI{}{\per\second}$  &$i$  & $a_i$ in $\SI{}{\per\second}$  \\ \hline
1    &   2556   &  2    &   2688   &  3    &   2786   &  4    &   2862   &  5    &   2944 \\ 
6    &   3015   &  7    &   3081   &  8    &   3142   &  9    &   3205   &  10    &   3243 \\ 
11    &   3290   &  12    &   3323   &  13    &   3365   &  14    &   3382   &  15    &   3375 \\ 
16    &   3387   &  17    &   3392   &  18    &   3367   &  19    &   3355   &  20    &   3342 \\ 
21    &   3330   &  22    &   3322   &  23    &   3314   &  24    &   3310   &  25    &   3307 \\ 
26    &   3312   &  27    &   3308   &  28    &   3311   &  29    &   3315   &  30    &   3334 \\ 
31    &   3327   &  32    &   3330   &  33    &   3332   &  34    &   3335   &  35    &   3338 \\ 
36    &   3345   &  37    &   3355   &  38    &   3350   &  39    &   3351   &  40    &   3352 \\ 
41    &   3354   &  42    &   3361   &  43    &   3360   &  44    &   3367   &  45    &   3367 \\ 
46    &   3368   &  47    &   3370   &  48    &   3367   &  49    &   3377   &  50    &   3371 \\ 
51    &   3373   &  52    &   3377   &  53    &   3379   &  54    &   3370   &  55    &   3372 \\ 
56    &   3375   &  57    &   3372   &  58    &   3377   &  59    &   3373   &  60    &   3373 \\ 
61    &   3373   &  62    &   3373   &  63    &   3385   &  64    &   3373   &  65    &   3390 \\ 
66    &   3375   &  67    &   3374   &  68    &   3376   &  69    &   3375   &  70    &   3378 \\ 
71    &   3375   &  72    &   3375   &  73    &   3377   &  74    &   3377   &  75    &   3382 \\ 
76    &   3375   &  77    &   3375   &  78    &   3377   &  79    &   3375   &  80    &   3383 \\ 
81    &   3377   &  82    &   3375   &  83    &   3377   &  84    &   3376   &  85    &   3382 \\ 
86    &   3377   &  87    &   3376   &  88    &   3377   &  89    &   3376   &  90    &   3383 \\ 
91    &   3378   &  92    &   3396   &  93    &   3376   &  94    &   3377   &  95    &   3376 \\ 
96    &   3377   &  97    &   3383   &  98    &   3378   &  99    &   3376   &  100    &   3377 \\ 
\end{tabular}\caption{\textbf{Fusion rate amplitudes.} Estimated peak heights of $f_{\text{Gaussians}}$. }\label{tab:ai}
\end{table}




\quad 

\noindent \textbf{Impulse response function.} 
The impulse response function was taken from Ref.~\cite{kobbersmed2020rapid}:
\begin{align}\label{eq:g}
    g(t) &= A\left(1-e^{-\frac{t-t_0}{\tau_{\text{r}}}}\right)\left(B e^{-\frac{t-t_0}{\tau_{\text{df}}}} + (1-B)e^{-\frac{t-t_0}{\tau_{\text{ds}}}}\right),
\end{align}
where $t_0=\SI{3}{\milli \second}$ is the onset, $A=\SI{7.21}{\micro \ampere}$ is the full amplitude (if there was no decay), $B=2.7\times 10^{-9}$ is the fraction of the fast decay, and $\tau_{\text{r}} = \SI{10.6928}{\second}$, $\tau_{\text{df}} = \SI{1.5}{\milli \second}$, $\tau_{\text{ds}} = \SI{2.8}{\milli \second}$ are the time constants of rise, fast decay and slow decay, respectively.

\subsection{Parameter studies}\label{sec:ParamStudies}
In order to contextualize our findings on the sensitivities as depicted in Fig.~\ref{fig:sens_normal}, we need to evaluate the range of possible system behaviors in response to different parameter values. For clarity, we limit our focus to alterations of the priming rate $k_R$, the value of which was previously found via optimization, and the two recovery rates  $g_V$ and $g_P$, which were estimated from the literature in our example (see Sec.~\ref{parameter_estimation}). 

\subsubsection{Varying the docking rate}

\begin{figure}[ht]
		\centering
		{\includegraphics[width=\textwidth]{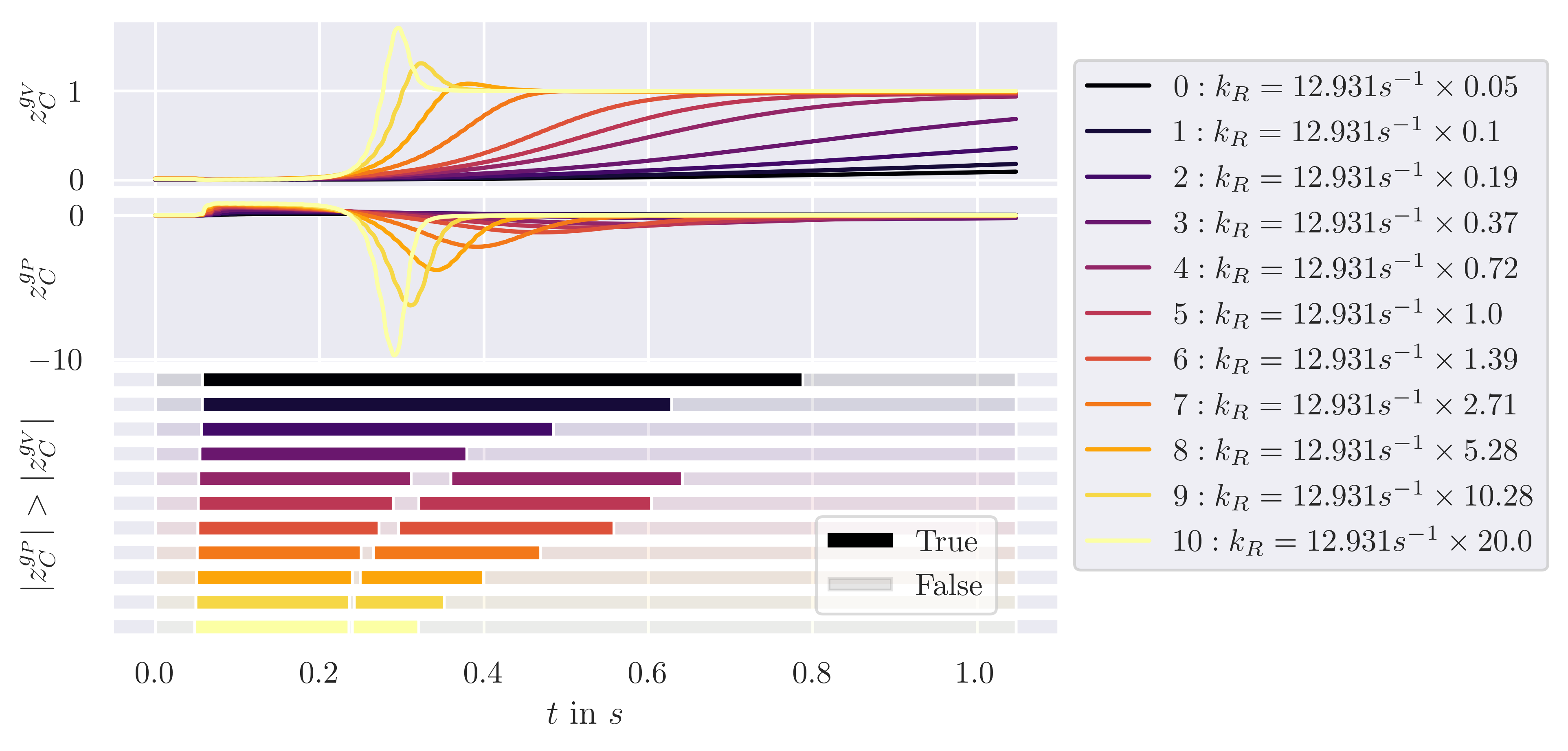}}
		\caption{\textbf{Impact of changing $k_R$ on the sensitivity time course.} The priming rate is varied from $1/20$ to $20$ times its original value while keeping all other parameter values as in our example. Note the logarithmic spacing between different values of $k_R$ (dark - low values, light - high values) and that the crimson color (line No. 5) corresponds to the parameter values used in the example in Fig.~\ref{fig:sens_normal}.}  
		\label{fig:kr_sweep}
	\end{figure}

The result of varying $k_R$ from $1/20$ to $20$ times its original value while keeping all other parameter values as in our example is depicted in Fig.~\ref{fig:kr_sweep}. Note the logarithmic spacing between different values of $k_R$ (dark - low values, light - high values) and that the crimson color (line No. 5) corresponds to the parameter values used in the example in Fig.~\ref{fig:sens_normal}. The top two graphs show the temporal evolution of the sensitivities, while the bar plots beneath give the behavior of the sensitivities' absolute values. At stimulation onset ($t_0=\SI{3}{\milli \second}$), the dominant sensitivity is $z^{g_P}_C$ for all values of $k_R$ under consideration, and after $\SI{1}{\second}$ of stimulation, $z^{g_P}_C$ always dominates, i.e. the identity of the rate-determining process is time-dependent for all examined values of $k_R$. 

Since the value of $k_R$ regulates the speed of the priming reaction for most of the stimulation time (as the unpriming rate $k_U$ falls to a negligible value after the first peak), one might naively expect a simple temporal compression (elongation) of the crimson system evolution in response to an increase (decrease) in $k_R$. While the sensitivity plots (top) do show this general behavior, interestingly, for high $k_R$, we also observe the formation of peaks of increased magnitude in the sensitivity graphs. The plots can be explained as follows: For small values of $k_R$, the priming reaction happens so slowly that neither release site nor vesicle recovery can develop much impact on the resulting weak signal within $\SI{1}{\second}$ and both sensitivities are generally small. Due to the availability of the vesicle reserve, release site recovery is the limiting process until $W_V$ has filled up sufficiently. At increased $k_R$, the vesicle supply $V$ is emptied at a greater speed that is dictated mainly by the recovery rate $g_p$. The amplitude of the resulting current $C$ is large as long as there are still vesicles available and exhibits a sharp decay after vesicle depletion -  the higher $k_R$ is, the steeper the decay. This is why a small increase in $g_P$ can lead to a strong relative attenuation of $C$, i.e. large negative peak values of $z^{g_P}_C$, at the end of vesicle depletion: increasing the recovery rate $g_P$ slightly shifts the time of vesicle depletion to the left, and the resulting relative difference in $C$ in the decay region is larger for a more steeply-decaying signal, leading to larger and sharper peaks for increased $k_R$. 

The peak emergence in $z^{g_V}_C$ can be explained in an analogous manner, however, since a slight increase in $g_V$ does not alter the vesicle depletion process much, the peaks are much smaller in magnitude (note the different scaling of the vertical axes). The combination of these effects leads to the formation of a second domain in which $z^{g_P}_C$ is the dominant sensitivity for a majority of the examined parameter space. After most vesicles have accumulated in $W_V$, the system is most sensitive to changes in $g_V$.

\subsubsection{Varying the vesicle recovery rate}

\begin{figure}[ht]
		\centering
		{\includegraphics[width=\textwidth]{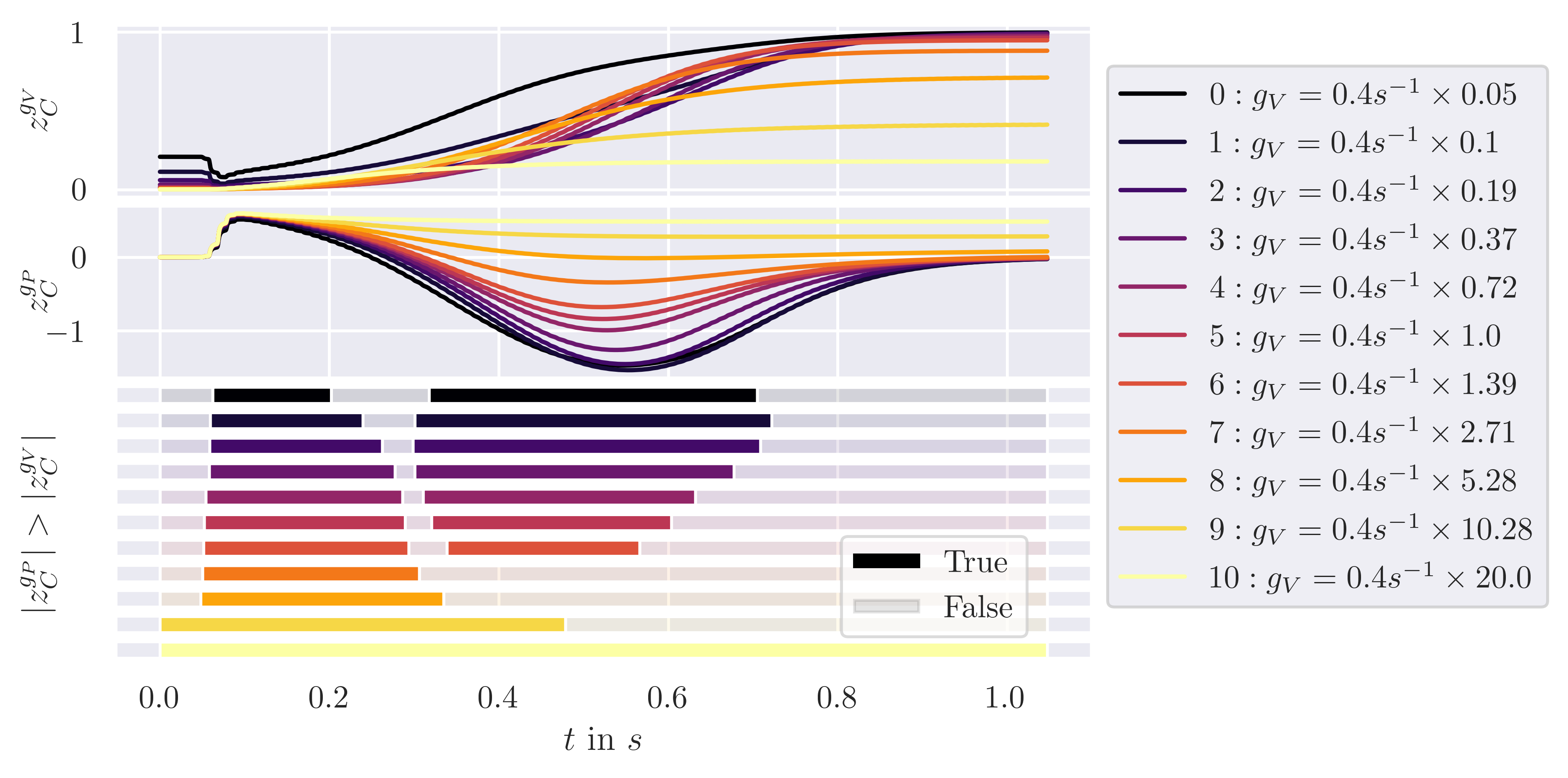}}
		\caption{\textbf{Impact of changing $g_V$ on the sensitivity time course.} The vesicle recovery rate is varied from $1/20$ to $20$ times its original value while keeping all other parameter values as in our example. Note the logarithmic spacing between different values of $g_V$ (dark - low values, light - high values) and that the crimson color (line No. 5) corresponds to the parameter values used in the example in Fig.~\ref{fig:sens_normal}.}  
		\label{fig:gv_sweep}
	\end{figure}

Fig.~\ref{fig:gv_sweep} shows the impact of of varying the vesicle recovery rate $g_V$ from $1/20$ to $20$ times its original value while keeping all other parameters as in the example. Increasing $g_V$ induces a flattening of the sensitivity $z^{g_V}_C$ curve to lower values, while the indentation in the course of $z^{g_P}_C$ gradually becomes less negative and finally levels out to an almost constant positive value. This is due to the fact that raising $g_V$ decreases vesicle depletion and thereby alleviates the effects discussed in the previous subsection. Vice versa, lowering $g_V$ increases vesicle depletion and its impact on the sensitivities. As a result, the identity of the rate-limiting process keeps behaving in a way similar to our example (again, lines/bar No. 5) for low values of $g_V$. Interestingly, when increasing $g_V$, the second domain of higher absolute sensitivity $z^{g_P}_C$ initially disappears and then, the first domain starts to expand significantly. Thus, the total amount of time spent in the site-limited state actually first decreases before increasing! This is especially relevant since it means that the same percentage of time spent in the site-limited state can result from different vesicle recovery rates. (If there was a way to distinguish the states experimentally, the short intermediate vesicle-limited domain may be too small to resolve.) Finally, when vesicles are replenished at very high speeds, the system is fully site-limited during the stimulation time. 

\subsubsection{Varying the release site recovery rate}

\begin{figure}[ht]
		\centering
        {\includegraphics[width=\textwidth]{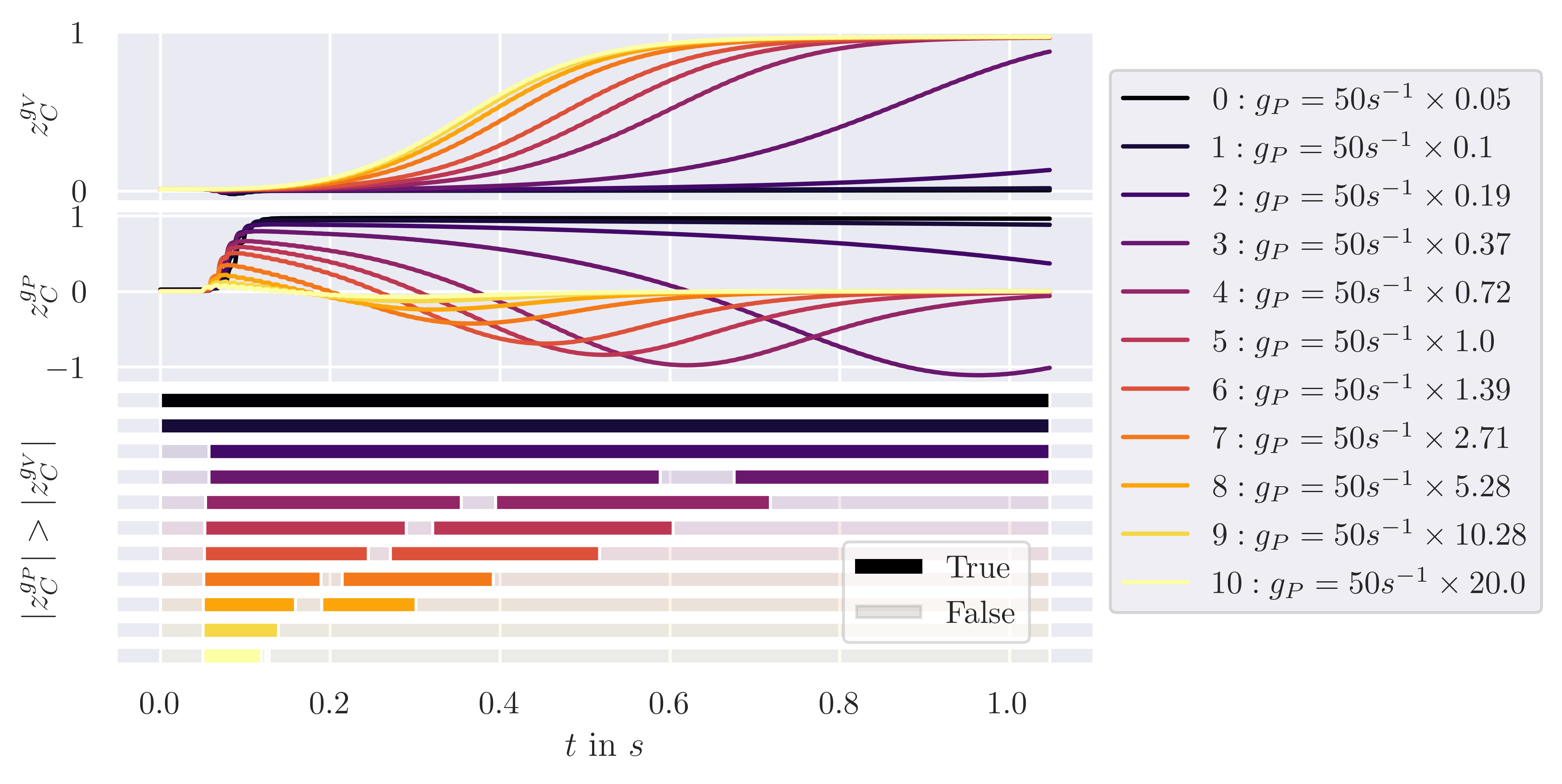}}
		\caption{\textbf{Impact of changing $g_P$ on the sensitivity time course.} The release site recovery rate is varied from $1/20$ to $20$ times its original value while keeping all other parameter values as in our example. Note the logarithmic spacing between different values of $g_P$ (dark - low values, light - high values) and that the crimson color (line No. 5) corresponds to the parameter values used in the example in Fig.~\ref{fig:sens_normal}.}  
		\label{fig:gP_sweep}
	\end{figure}

The effect of varying the release site recovery rate $g_P$ from $1/20$ to $20$ times its original value while keeping all other parameters as in the example is depicted in Fig.~\ref{fig:gP_sweep}. An increase in $g_P$ results in temporal compression of both sensitivities. For the sensitivity $z^{g_V}_C$, a decrease in $g_P$ simply has the opposite effect of a temporal stretching of the time course. This is because the release site recovery rate determines how fast the vesicle supply is emptied and $W_V$ is filled, and the earlier this happens, the earlier the sensitivity to vesicle recovery rises. For $z^{g_P}_C$, raising $g_P$ also brings on a strong amplitude diminution while lowering $g_P$ has the opposite effect. At high release site recovery rates, vesicles are depleted quickly and the resulting signal $C$ shows an exponential decay. A small increase in $g_p$ steepens the slope of this decay, however, there is a limit to this steepening since the vesicle depletion process is still constrained by the amount of priming and fusions that happen. Thus, at high $g_P$, the slope changes only very slightly which is why the relative change in the current $C$ and therefore also the sensitivity $z^{g_P}_C$ is small. At low values of $g_P$, as vesicle depletion happens very slowly, site recovery speed has the greatest impact on signal strength and even small increases in $g_P$ can result in a lasting stronger signal. 

The combination of these effects results in the behavior of the limiting process that is depicted on the bottom of Fig.~\ref{fig:gP_sweep}: Except for very low values of $g_P$, the identity switches at least once and always begins as site-limiting at stimulation onset, before $W_V$ has filled sufficiently for vesicle recovery to have an impact. The two domains from our example where the system is site-limited are conserved within a range of $g_P$ but are compressed with increasing $g_P$. For very high $g_P$, the second domain disappears completely and the system is only site-limited for a short amount of time at stimulation onset. Only at very low site recovery rates, site recovery is always the rate-determining process.

\subsection{Steady state investigation} \label{app:steady_state}

In the following we show that for constant $k_F> 0$ and constant $k_U >0$ the system given by Eq.~\eqref{RRE} has a unique steady state. Assume also $g_V,g_P>0$. 
The steady state $ \hat{\boldsymbol{x}}=(\hat{V},\hat{W}_V,\hat{W}_P,\hat{R},\hat{P})$ is given by the fixed point equations
\begin{align}
    0 & = -k_R\hat{V}\hat{P} + g_V \hat{W}_V + k_U \hat{R} \label{steady_V} \\
    0 & = k_F \hat{R} - g_V \hat{W}_V \label{steady_WV}  \\
    0 & = k_F \hat{R} - g_P \hat{W}_P \label{steady_WP} \\
    0 & = k_R \hat{V}\hat{P} - k_U \hat{R} -k_F \hat{R} \label{steady_R} \\
    0 & = -k_R \hat{V}\hat{P} + g_P \hat{W}_P +k_U\hat{R} \label{steady_P} 
\end{align}
with 
\begin{align}
    \hat{R} + \hat{P} +\hat{W}_P& = n_{\text{sites}}, \label{sum_sites} \\
    \hat{V} + \hat{P} + \hat{W}_V & = n_{\text{ves}}. \label{sum_ves}
\end{align}
Note that \eqref{steady_R} follows from \eqref{steady_V} and \eqref{steady_WV}, while \eqref{steady_P} follows from \eqref{steady_WP} and \eqref{steady_R}, so both \eqref{steady_R} and \eqref{steady_P} are redundant. 

From \eqref{steady_WV} it follows $\hat{W}_V=\frac{k_F}{g_V}\hat{R}$ and from  \eqref{steady_WP} it follows $\hat{W}_P=\frac{k_F}{g_P}\hat{R}$. Inserting into \eqref{sum_sites} and \eqref{sum_ves} we get
\begin{equation}\label{eq:P2_SS}
    \hat{P} =  n_{\text{sites}} - \left(1+\frac{k_F}{g_P}\right)\hat{R}
\end{equation}
and 
\begin{equation}
    \hat{V} = n_{\text{ves}} - \left(1+\frac{k_F}{g_V}\right)\hat{R},
\end{equation}
respectively. Set $\alpha:=1+\frac{k_F}{g_P}> 1 $, $\beta := 1+\frac{k_F}{g_V} > 1$ and $\gamma:=\frac{k_F+k_U}{k_R}> 0$. 
Inserting into \eqref{steady_V} we obtain
\begin{equation}
    0 = \hat{R}^2  \underbrace{-\left(\frac{n_{\text{ves}}}{\beta}+ \frac{n_{\text{sites}}}{\alpha} + \frac{\gamma}{\alpha \beta}\right)}_{=:p} \hat{R} + \underbrace{\frac{n_{\text{ves}}n_{\text{sites}}}{\alpha \beta}}_{=:q}.
\end{equation}
This yields two solutions for $\hat{R}$:
\begin{equation}
            \hat{R}_{1,2} = \underbrace{-\frac{p}{2}}_\text{=:$R_L$} \pm \underbrace{\sqrt{\frac{p^2}{4}-q}}_\text{=:$R_R$} = R_L \pm R_R.
\end{equation}
We will now show that one of these solutions can be discarded as it leads to negative values for $\hat{P}$. From \eqref{eq:P2_SS} we have
        \begin{align*}
            \hat{P}_{1,2}&= n_{\text{sites}} - \alpha \hat{R}_{1,2} =  n_{\text{sites}} - \alpha (R_L \pm R_R) =   n_{\text{sites}} - \alpha R_L \mp \alpha R_R.
        \end{align*}
For $\alpha R_R > n_{\text{sites}} - \alpha R_L$ this will give $\hat{P}_1<0$. Let us therefore compare the two expressions in the following. It holds
\begin{align}
    \lefteqn{n_{\text{sites}} - \alpha R_L } \nonumber\\
    & = n_{\text{sites}} + \alpha \frac{p}{2}  \nonumber \\
    & = \frac{\beta  n_{\text{sites}} - \alpha  n_{\text{ves}}-\gamma}{2 \beta } \nonumber \\
    & = \sqrt{\frac{(\beta  n_{\text{sites}} - \alpha  n_{\text{ves}}-\gamma)^2}{(2 \beta )^2}} \nonumber\\
    & = \sqrt{\frac{(\beta  n_{\text{sites}})^2 - 2\beta  n_{\text{sites}}( \alpha  n_{\text{ves}}+\gamma) + ( \alpha  n_{\text{ves}}+\gamma)^2}{(2 \beta )^2}} \label{left}
\end{align}
and 
\begin{align}
   \lefteqn{ \alpha R_R} \nonumber\\
  & = \sqrt{\alpha^2 \left(\frac{p^2}{4} -q\right)} \nonumber\\
     & = \sqrt{\frac{(\beta  n_{\text{sites}} + \alpha  n_{\text{ves}}+\gamma)^2 - 4 \alpha \beta  n_{\text{ves}} n_{\text{sites}}}{(2\beta )^2}} \nonumber\\
         & = \sqrt{\frac{(\beta n_{\text{sites}})^2 + 2\beta n_{\text{sites}}( \alpha  n_{\text{ves}}+\gamma) + ( \alpha  n_{\text{ves}}+\gamma)^2 - 4 \alpha \beta  n_{\text{ves}} n_{\text{sites}}}{(2\beta )^2}} \nonumber\\
                  & = \sqrt{\frac{(\beta  n_{\text{sites}})^2 - 2\beta n_{\text{sites}}( \alpha  n_{\text{ves}}-\gamma) + ( \alpha  n_{\text{ves}}+\gamma)^2 }{(2\beta )^2}}. \label{right}
\end{align}
Since $\gamma >0 $, comparing \eqref{left} to \eqref{right} proves that indeed $\alpha R_R > n_{\text{sites}} - \alpha R_L$ and we need to choose $\hat{R}_2$: 
\begin{align}
 \hat{R} & = \hat{R}_2 =  -\frac{p}{2} - \sqrt{\frac{p^2}{4}-q} \nonumber\\
    & = \frac{\beta  n_{\text{sites}}+ \alpha  n_{\text{ves}}+\gamma}{2 \alpha\beta } - \sqrt{\frac{(\beta  n_{\text{sites}}+ \alpha  n_{\text{ves}}+\gamma)^2}{(2 \alpha\beta )^2} - \frac{n_{\text{ves}}n_{\text{sites}}}{\alpha \beta}}.
\end{align}
We note that the term under the square root is always non-negative since
\begin{align}
    & \frac{(\beta  n_{\text{sites}}+ \alpha n_{\text{ves}}+\gamma)^2}{(2 \alpha\beta)^2} - \frac{n_{\text{ves}}n_{\text{sites}}}{\alpha \beta}  \nonumber \\
    & = \frac{1}{(2 \alpha\beta )^2} \left[(\beta  n_{\text{sites}}+ \alpha n_{\text{ves}})^2 + 2 (\beta  n_{\text{sites}}+ \alpha  n_{\text{ves}})\gamma + \gamma^2 - 2 \alpha \beta n_{\text{ves}}n_{\text{sites}} \right] \nonumber\\
      & = \frac{1}{(2 \alpha\beta )^2} \left[(\beta  n_{\text{sites}})^2+ (\alpha n_{\text{ves}})^2 + 2 (\beta n_{\text{sites}}+ \alpha  n_{\text{ves}})\gamma + \gamma^2 \right]  \nonumber \\
       & \geq 0.
\end{align}
In summary this means that for each choice of (positive) parameter values there exist two fixed points, but only one with physically relevant numbers while the other one has negative values. That is, there is a unique steady state of the system. Due to the stoichiometric structure of the system, this steady state will be approached in the course of time, no matter which initial state (with non-negative values) is chosen.

For time-dependent, periodic rates $k_F(t)$ the system will be pulled towards the time-dependent steady state, thereby showing itself a periodic behavior, see Fig.~\ref{fig:normal}.

\subsection{Reduced reaction systems}
\label{reduced_reaction_system}

In Sec.~\ref{sec:moments} we have argued that the similarity of the ODE-solution with the stochastic mean stems from the independence of the recovery steps. We now clarify this issue by comparing two reduced reaction systems:
\begin{enumerate}
    \item[(I)]  Standard binding and unbinding given by the reactions
\begin{equation}
    A+B \stackrel{\alpha}{\longrightarrow} W, \quad W \stackrel{\beta}{\longrightarrow} A+B
    \label{eq:systemI}
\end{equation}
with the associated ODEs given by 
\begin{equation}
    \frac{d}{dt}A(t)=\frac{d}{dt}B(t)= -\alpha A(t)B(t)+\beta W(t) = -\frac{d}{dt}W(t).
\end{equation}
\item[(II)] Binding with independent return given by the reactions
\begin{equation}
     A+B \stackrel{\alpha'}{\longrightarrow} W_A +W_B, \quad W_A \stackrel{g_A}{\longrightarrow} A, \quad W_B \stackrel{g_B}{\longrightarrow} B
     \label{eq:systemII}
\end{equation}
with the associated ODEs given by
\begin{align}
    \frac{d}{dt}A(t) & = -\alpha' A(t)B(t)+g_A W_A(t) = -\frac{d}{dt}W_A(t), \\
    \frac{d}{dt}B(t) &= -\alpha' A(t)B(t)+g_B W_B(t) = -\frac{d}{dt}W_B(t).
\end{align}
\end{enumerate}
We note that system (I) arises from our full reaction system (as depicted in Fig.~\ref{fig:simplified_model}) by setting $g_V=g_P=\infty$, $k_U=0$ $k_F=\beta$, $k_R=\alpha$ (with the species being related by $V=A$, $P = B$, $R= W$), while the second system (II) results from setting $k_F=\infty$ and $k_U=0$, $k_R=\alpha'$, $g_V=g_A$, $g_P=g_B$ (with the species being related by $V=A$, $P=B$, $W_V=W_A$, $W_P=W_B$). 

\begin{figure}[ht]
			\centering
	\begin{subfigure}[b]{0.49\textwidth}
	    	\centering
	\includegraphics[width=1.\textwidth]{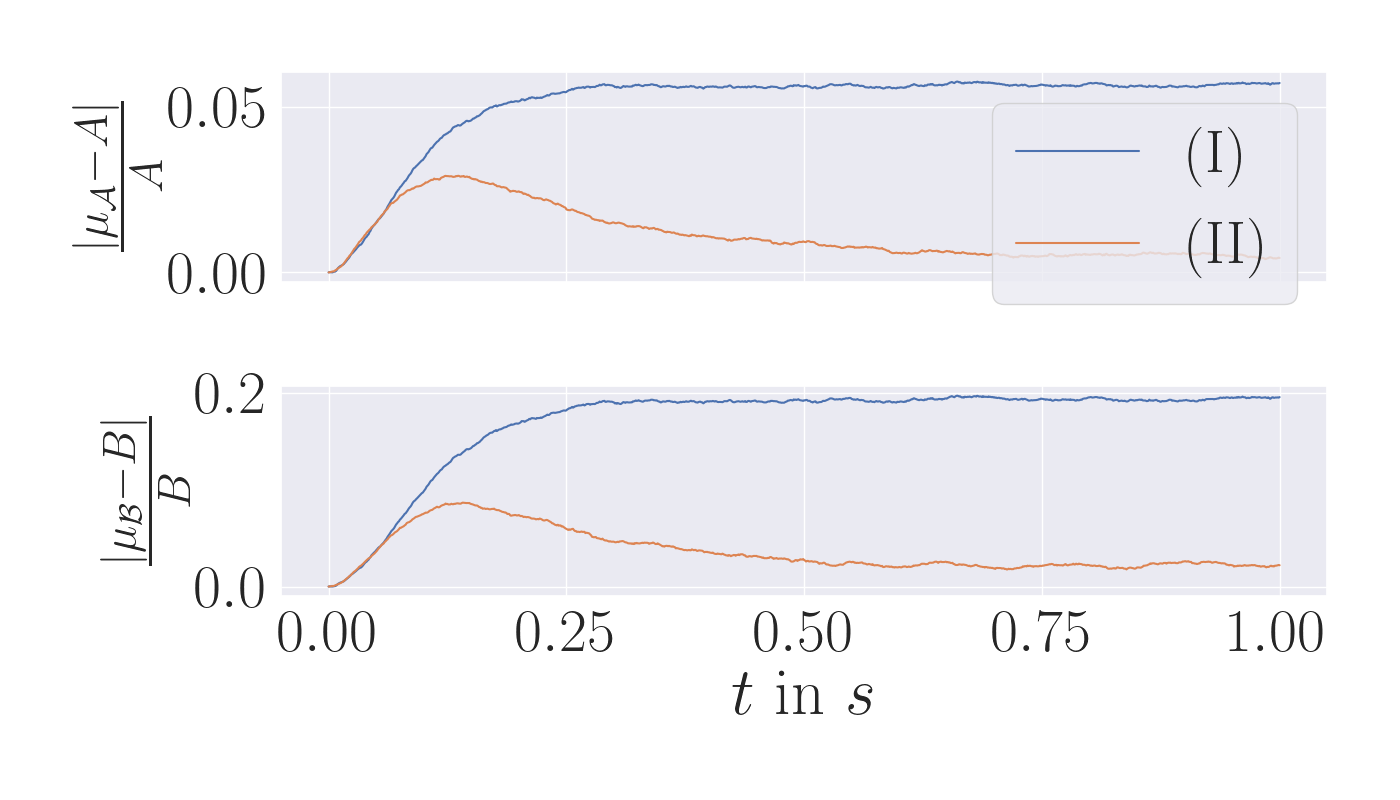}
	\caption{Relative error}\label{fig:rel_error_two_systems}
	\end{subfigure}
	\hfill
		\begin{subfigure}[b]{0.49\textwidth}
	    	\centering
	\includegraphics[width=1.\textwidth]{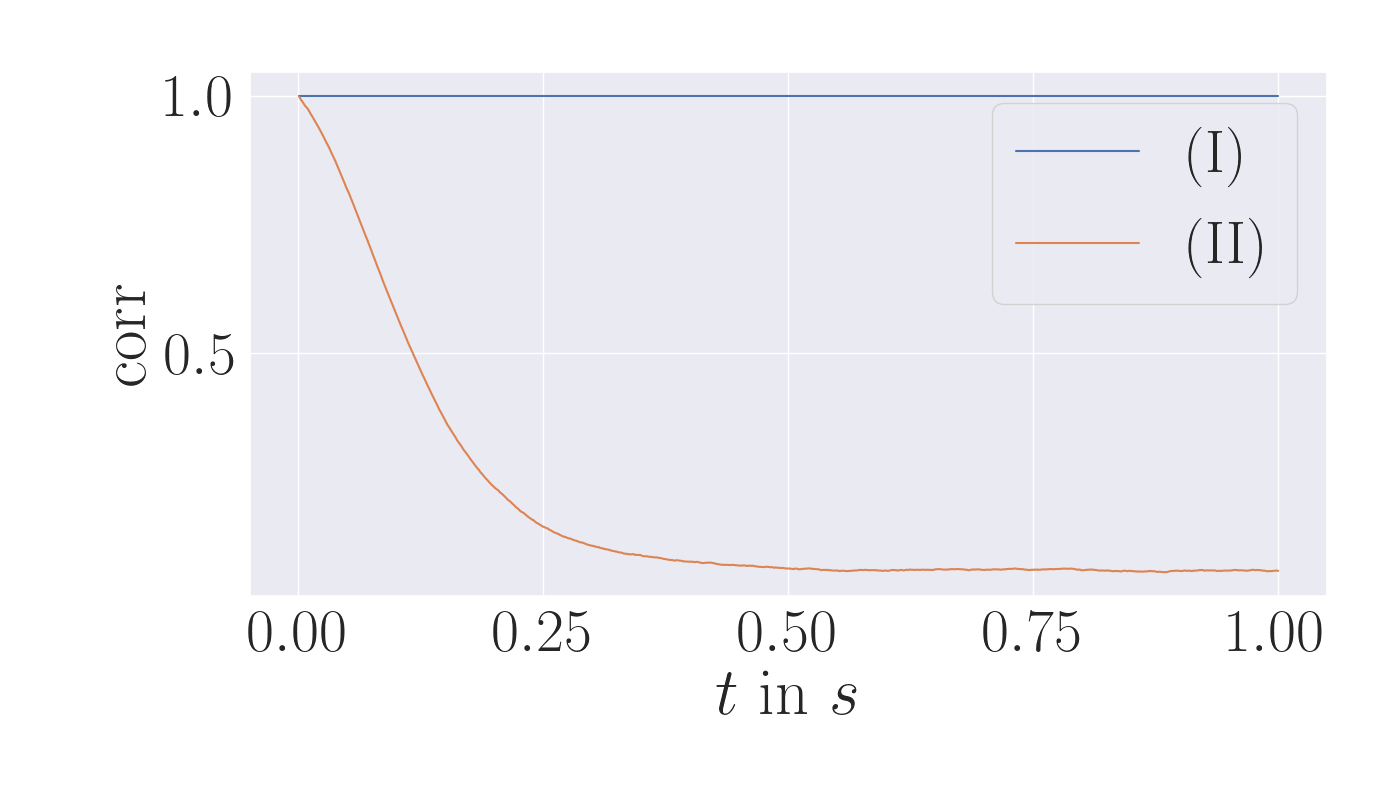}
	\caption{Correlation} \label{fig:correlation_two_systems}
	\end{subfigure}
		\caption{Comparison of reduced artificial systems (I) and (II) from Equations \eqref{eq:systemI} and \eqref{eq:systemII}. (a) Relative error $\vert\mu_{\mathcal{X}_i}(t)-X_i(t)\vert/X_i(t)$ for $\mathcal{X}_i=\mathcal{A}, \mathcal{B}$ estimated from $2\cdot 10^5$ simulations. (b) Correlation function $\text{corr}(t)$ defined in \eqref{corr_AB}. 
        Paramter values: 
      System (I):   $\alpha= \beta=5$. System (II): $\alpha'=g_A=g_B=5$. Both systems with initial states $A(0)=2, B(0)=1$ and all other species start in $0$. 
  }
	\end{figure}

It can easily be shown that for $\alpha=\alpha'$, $\beta=g_A=g_B$ and appropriate initial states (satisfying $W(0)=W_A(0)=W_B(0)$), the ODE-solutions of the two systems (I) and (II) fully agree. However, the first-order moments of the corresponding stochastic jump processes $\mathcal{A}(t), \mathcal{B}(t)$ are not the same. Actually, the first-order moments $\mu_{\mathcal{A}}(t), \mu_{\mathcal{B}}(t)$ of the second system (II) of independent return are closer to the ODE-solution $A(t), B(t)$, see Fig.~\ref{fig:rel_error_two_systems}, where we plot the relative errors. Denote the corresponding standard deviations of $\mathcal{A}(t)$ and $\mathcal{B}(t)$ by $\sigma_\mathcal{A}(t), \sigma_\mathcal{B}(t)$, respectively. Then Fig.~\ref{fig:correlation_two_systems} shows the correlation function
\begin{equation}\label{corr_AB}
    \text{corr}(t) := \frac{\text{cov}(\mathcal{A}(t),\mathcal{B}(t))}{\sigma_\mathcal{A}(t)\sigma_\mathcal{B}(t)}
\end{equation}
for both reduced systems, with significantly smaller values for the second system (II). This confirms our hypothesis that independent recovery increases the approximation quality of the ODE-system to the mean of the stochastic dynamics. 
	
\bibliography{bibliography}
\bibliographystyle{unsrt}

\end{document}


\begin{table}[ht]
\begin{tabular}{l l||l l||l l||l l||l l}
$i$ & $a_i$ in $\SI{}{\per\second}$  & $i$  & $a_i$ in $\SI{}{\per\second}$ &$i$  & $a_i$ in $\SI{}{\per\second}$  &$i$  &$a_i$ in $\SI{}{\per\second}$  &$i$  & $a_i$ in $\SI{}{\per\second}$  \\ \hline
1    &   2556.0   &  2    &   2688.0   &  3    &   2786.0   &  4    &   2862.0   &  5    &   2944.0 \\ 
11    &   3290.0   &  12    &   3323.0   &  13    &   3365.0   &  14    &   3382.0   &  15    &   3375.0 \\ 
21    &   3330.0   &  22    &   3322.0   &  23    &   3314.0   &  24    &   3310.0   &  25    &   3307.0 \\ 
31    &   3327.0   &  32    &   3330.0   &  33    &   3332.0   &  34    &   3335.0   &  35    &   3338.0 \\ 
41    &   3354.0   &  42    &   3361.0   &  43    &   3360.0   &  44    &   3367.0   &  45    &   3367.0 \\ 
51    &   3373.0   &  52    &   3377.0   &  53    &   3379.0   &  54    &   3370.0   &  55    &   3372.0 \\ 
61    &   3373.0   &  62    &   3373.0   &  63    &   3385.0   &  64    &   3373.0   &  65    &   3390.0 \\ 
71    &   3375.0   &  72    &   3375.0   &  73    &   3377.0   &  74    &   3377.0   &  75    &   3382.0 \\ 
81    &   3377.0   &  82    &   3375.0   &  83    &   3377.0   &  84    &   3376.0   &  85    &   3382.0 \\ 
91    &   3378.0   &  92    &   3396.0   &  93    &   3376.0   &  94    &   3377.0   &  95    &   3376.0 \\ 
96 & 3377.0 & 97 & 3383.0 & 98 & 3378.0 & 99 & 3376.0 & 100 &
3377.0
\end{tabular}
\end{table}